%% file: arxiv_main.tex
\documentclass[final]{arxiv} 

\OneAndAHalfSpacedXI 

\usepackage[utf8]{inputenc} 
\usepackage[T1]{fontenc}    
\usepackage{bm}
\usepackage{nicefrac}       
\usepackage{microtype}      
\usepackage{xcolor}         
\usepackage{colortbl}       
\usepackage{subfigure}
\usepackage{booktabs}       
\usepackage{multirow}
\usepackage{tablefootnote}
\usepackage{float}
\usepackage{algorithm}
\usepackage{algpseudocode}
\usepackage{tikz}
\usepackage{pict2e}
\usepackage{hyperref}       
\usepackage{xr-hyper}       
\usepackage{cleveref}       

\algnewcommand\algorithmicswitch{\textbf{switch}}
\algnewcommand\algorithmiccase{\textbf{case}}
\algnewcommand\algorithmicdefault{\textbf{default}}
\algdef{SE}[SWITCH]{Switch}{EndSwitch}[1]{\algorithmicswitch\ #1\ \algorithmicdo}{\algorithmicend\ \algorithmicswitch}
\algdef{SE}[CASE]{Case}{EndCase}[1]{\algorithmiccase\ #1:}{}
\algtext*{EndCase}
\algtext*{EndSwitch}

\algrenewcommand{\algorithmiccomment}[1]{\hfill\textcolor{gray}{\textit{// #1}}}

\newcommand\Set[1]{\ensuremath{\mathcal{#1}}}

\newcommand{\SetJ}{\Set{J}}
\newcommand{\SetT}{\Set{T}}
\newcommand{\SetK}{\Set{K}}
\newcommand{\SetS}{\Set{S}}

\usepackage{natbib}
\bibpunct[, ]{(}{)}{,}{a}{}{,}%
\TheoremsNumberedThrough     
\ECRepeatTheorems

\EquationsNumberedThrough    

\begin{document}


\RUNAUTHOR{Maragno et al.}

\RUNTITLE{Learning-Augmented Optimization for Strategic Two-Echelon Spare Parts Network Design}

\TITLE{Learning-Augmented Optimization for Strategic Two-Echelon Spare Parts Network Design}
\ARTICLEAUTHORS{%
\AUTHOR{Donato Maragno,\textsuperscript{a} Marco Caserta,\textsuperscript{a,b} Alberto Sinigaglia,\textsuperscript{c} Komlanvi Ametana,\textsuperscript{d} David Corredor Montenegro,\textsuperscript{a} Luca D'Angelo\textsuperscript{a}} \AFF{\textsuperscript{a}Amazon RME, \EMAIL{dmaragno@amazon.de, mcaserta@amazon.lu,
davcrd@amazon.com, dangellu@amazon.lu}} 
\AFF{\textsuperscript{b}IE University, Department of Operations, \EMAIL{marco.caserta@ie.edu}}
\AFF{\textsuperscript{c}Reinforcement Learning Research Lab, UniPD, \EMAIL{alberto.sinigaglia@phd.unipd.it}} 
\AFF{\textsuperscript{d}\EMAIL{kametanap@gmail.com}}}

\ABSTRACT{%

  We study the strategic design of a two echelon spare parts
  inventory network in which evaluating each candidate topology
  requires an expensive inventory optimization model. The design
  problem partitions hundreds of sites into feasible clusters and
  selects a central replenishment site for each cluster to reduce
  purchasing and inventory costs while maintaining service levels.
  This creates a broader challenge for learning-augmented
  optimization: because the optimizer preferentially explores
  candidates with high predicted savings, it can systematically
  exploit optimistic errors in surrogate predictions.

  We develop a conservative framework combining a graph neural
  network ensemble, variable neighborhood search, and set
  partitioning recombination. The surrogate is trained on exact
  cluster evaluations, while a lower quantile of the economic savings
  predicted by the ensemble guides the search to reduce the influence
  of optimistic predictions. Clusters found during the search are
  recombined through a set partitioning model with objective
  coefficients given by the surrogate. The resulting network is
  evaluated with the exact inventory model, and only this evaluation
  is used to report performance.

  In a case study of 246 fulfillment centers in Amazon’s North
  American network, the full framework improves combined savings by
  30.5\% relative to an optimization baseline constructed entirely
  from exact cluster evaluations. Across six independent
  replications, the framework achieves this improvement while
  maintaining a service level of approximately 99.8\%. Component
  experiments under equal computational budgets provide descriptive
  evidence that the search guided by the graph surrogate produces
  higher mean exact savings than the search guided by a tabular
  alternative under both mean and conservative scoring. Conservative
  scoring also improves mean savings for both surrogate classes and
  reduces the share of clusters selected into the final network that
  are overestimated by the graph surrogate from 68\% to 28\%.
  Predictive and ranking accuracy deteriorate among
  search-generated candidates that receive high surrogate scores,
  indicating that performance on a random holdout sample can
  provide an incomplete picture of surrogate performance during
  optimization. The results indicate that graph representation
  improves the ranking of candidates relevant to the decision,
  conservative scoring limits optimistic prediction errors, and
  exact validation ensures that reported performance does not
  depend on surrogate predictions.
}%




\KEYWORDS{Spare parts inventory, network design, learning-augmented
optimization, surrogate-assisted optimization, graph neural networks}

\maketitle

\section{Introduction}
\label{sec:introduction}

Managing spare parts in large operational networks is a difficult and
strategically important inventory problem. Unlike conventional
consumer-goods settings, spare parts are characterized by highly
intermittent demand, long and uncertain replenishment lead times,
heterogeneous criticality, and often high unit costs
\citep{kennedy_overview_2002,Topan2020review}. Stockouts can delay
repairs and reduce equipment availability, while excess inventory
ties up substantial capital \citep{muckstadt2005analysis}. Firms
operating large installed equipment therefore face a persistent
tradeoff between maintaining high service levels and controlling
inventory costs. This challenge is particularly acute in
geographically distributed networks, where thousands of low-demand
but operationally critical items must be positioned across many sites
\citep{tapia-ubeda_modelling_2020}.

This paper is motivated by the redesign of Amazon's spare parts
inventory network. Amazon manages spare parts in more than 1{,}400
facilities, including fulfillment centers, sort centers, and delivery
stations, and its portfolio comprises more than 1.2 million parts
with annual purchase-order spend exceeding \$1 billion. In the current
system, sites operate largely as independent single-echelon entities
that replenish directly from external suppliers under min/max
inventory policies, also known as $(s,S)$ policies
\citep{zheng_finding_1991,axsater_inventory_2015}. Amazon is moving
toward a two-echelon structure in which selected central sites
replenish downstream leaf sites. If well designed, this transition
can shorten effective replenishment lead times, reduce inventory and
procurement costs, and maintain high service levels
\citep{sherbrooke_metric_1968,Howard2015}. The key strategic question
is which sites should serve as central sites and how the remaining
sites should be assigned to them.

Answering this question is difficult because the operational and
economic value of a candidate
topology cannot be assessed on geographical grounds alone. Strategic
clustering decisions determine lower-level inventory policies,
replenishment quantities, and shipment flows that ultimately drive
purchase-order spend, inventory levels, and service performance.
Although these decisions could in principle be modeled jointly
\citep{shen2003joint,tapia-ubeda_modelling_2020}, doing so is
computationally prohibitive at the scale considered here.
Learning-augmented optimization offers a potential way to overcome
this difficulty by using surrogate models when exact evaluations are
expensive \citep{elmachtoub2022smart,ChanLinSaxe2025}. Yet there is
limited work on large-scale strategic spare parts network design in
which candidate topologies must be evaluated through a
lower-level inventory model. Among the closest related works,
\citet{SpieckermannMinnerSchiffer2025} embed graph neural networks in
an optimization pipeline for a network design problem, but their
setting does not involve spare parts or expensive lower-level
inventory evaluation.

Using a surrogate inside an optimization procedure creates an
additional challenge: the search does not evaluate candidates at
random. It favors configurations with high predicted savings and can
therefore amplify optimistic prediction errors into poor design
decisions. This raises three questions: how can learned guidance be
combined with exact optimization to make large-scale search
computationally practical
without sacrificing decision quality; how should candidate network
configurations be represented by the surrogate; and how should
surrogate predictions be translated into search decisions when
optimization can amplify favorable prediction errors?

We address these questions with a learning-augmented framework for
large-scale strategic two-echelon spare parts network design.
Candidate clusters are
represented as attributed graphs and evaluated by an ensemble of GNN
surrogates trained on exact lower-level evaluations. A variable
neighborhood search uses these surrogates to explore a much larger
design space, scoring candidates with a conservative lower quantile
of the ensemble predictions. The clusters discovered during the
search are then recombined through a set partitioning model solved to optimality
and re-evaluated using the exact lower-level inventory model.
Learning therefore guides exploration and recombination, while exact
lower-level re-evaluation validates the reported performance of the
selected network.

This paper makes three contributions. First, we develop a
learning-augmented optimization approach for network design problems
with expensive lower-level evaluation. The approach uses exact
evaluations to train a surrogate, surrogate guidance to explore a
much larger design space, and exact lower-level evaluation to assess
the final design.

Second, we evaluate how surrogate choice affects downstream decision
quality in combinatorial network design. The GNN uses the graph
structure of candidate clusters to learn their economic value, and
our component analysis evaluates it not only through predictive
accuracy but through the exact quality of the designs produced during search.

Third, we examine how conservative surrogate scoring affects
prescriptive performance. Because the search preferentially explores
candidates with high predicted values, optimistic prediction errors
are more likely to influence the search. We therefore compare
candidates using a lower quantile of the ensemble predictions rather
than the mean.

We validate the framework in an industrial-scale case study based on
Amazon’s North American fulfillment-center network. Across six
independent replications, the full framework improves combined
savings by 30.5\% relative to an optimization-based baseline
constructed from exact evaluations while maintaining an average
service level of 99.8\%. Equal-budget component experiments assess
the effects of graph-based representation and conservative
surrogate scoring on the exact quality of the network designs
produced by the search. We also evaluate surrogate performance
under geographic distribution shift and on candidates generated by
the search, providing evidence on how predictive performance
changes in the regions of the design space most relevant to optimization.

The paper has the following structure.
Section~\ref{sec:literature-review} reviews the relevant literature
on multi-echelon spare parts systems and learning-augmented
optimization. Section~\ref{sec:problem-settings} defines the problem
setting and decision structure. Section~\ref{sec:nodal-model}
describes the lower-level two-echelon evaluation model.
Section~\ref{sec:solution-framework} presents the learning-augmented
solution framework. Section~\ref{sec:computational}
reports the computational study and evaluates predictive performance,
generalization, and downstream decision quality.
Section~\ref{sec:conclusions} concludes with managerial implications
and directions for future research.

\section{Literature Review}
\label{sec:literature-review}

Our paper lies at the intersection of two literature streams:
multi-echelon spare parts inventory systems and learning-augmented
optimization for problems with expensive exact evaluation. We review
these streams in turn and position our contribution at their intersection.

\subsection{Multi-echelon spare parts inventory systems}

The multi-echelon inventory literature starts with the seminal work
of \citet{ClarkScarf1960}, who establish the structural foundations
for analyzing inventory decisions across multiple stages. In spare
parts settings, these ideas were specialized and extended by
\citet{sherbrooke_metric_1968} and later synthesized in the
service-parts literature by \citet{muckstadt2005analysis}. Reviews by
\citet{kennedy_overview_2002}, \citet{bacchetti_spare_2012},
\citet{Eruguz2016}, \citet{Topan2020review},
\citet{zhang_spare_2021}, and \citet{caserta:dangelo:25} show that
spare parts systems differ from standard retail settings because
demand is highly intermittent, item criticality is heterogeneous, and
the economic consequences of stockouts are often severe.

A recurring theme in this literature is the value of centralization
and risk pooling. \citet{eppen_noteeffects_1979} shows how demand
aggregation can reduce expected mismatch costs, while related work
studies lateral transshipments and emergency supply options for
improving service performance under decentralized inventories
\citep{grahovac_sharing_2001,axsater_new_2003,paterson_inventory_2011,wong_multi-item_2006}.
In spare parts environments, however, such flexibility can increase
operational complexity. This motivates two-echelon structures in
which a central site supplies a set of downstream sites and captures
part of the pooling benefit through a more coordinated architecture.

Several papers analyze two-echelon spare parts systems from an
operational or tactical perspective. \citet{Dada1992} studies
priority shipments; \citet{Howard2015} examine emergency stocks and
pipeline information; and \citet{DrentArts2021} study expediting
decisions. Other work considers multi-item control under service
constraints \citep{Topan2010,Topan2017}, proactive and reactive
interventions \citep{TopanVanDerHeijden2020}, and order allocation,
inventory placement, stock allocation, and storage configuration
\citep{Wang2019,Wang2022,ChenDistribution2022,Dui2023}. Broader
frameworks and reviews are provided by \citet{Driessen2015} and
\citet{Topan2020review}, while recent studies explore reinforcement
learning for multi-echelon inventory control
\citep{Geevers2024,Stranieri2024}. This literature develops
increasingly sophisticated methods for operating multi-echelon
systems, but largely treats the network structure as given.

A related line of work considers supply chain design and safety-stock
placement. \citet{GravesWillems2000,GravesWillems2003} show how
strategic inventory positioning can be optimized over a supply
network, \citet{hammami:frein:14} study capacitated multi-echelon
inventory placement under lead-time constraints, and
\citet{shen2003joint} integrate location and inventory
decisions in a joint framework. For spare parts specifically,
\citet{tapia-ubeda_modelling_2020} study supply chain network design
problems that combine structural and inventory considerations. These
papers move beyond pure control decisions and recognize the
importance of network structure, but they do not address our central
computational setting: a large combinatorial network design problem
in which each candidate network must be evaluated through an
expensive lower-level two-echelon optimization model. We build on the
tactical and operational model developed in \citet{Caserta2026} but
make network topology, rather than inventory control within a given
topology, the primary decision.

\subsection{Learning-augmented optimization for expensive network design}

A separate but increasingly relevant stream studies how machine
learning and data analytics can support operational decision making
and optimization
\citep{feng:shanthikumar:23}. Within this literature,
\citet{elmachtoub2022smart} provide a conceptual
foundation by emphasizing that predictive models embedded in
optimization should be evaluated with downstream decision quality in
mind rather than by prediction error alone, while \citet{biggs+:23}
study constrained optimization directly over objective functions
represented by trained random forests.

Related work considers bilevel and hierarchical optimization, in which
upper-level design decisions interact with computationally expensive
lower-level models. Exact reformulations and decomposition methods
can become difficult to scale when the lower-level problem must be
solved repeatedly, motivating approximation and learning-based
approaches
\citep{sinha2016optimistic,sinha2018kriging,dumouchelle2024neurbilo,ChanLinSaxe2025}.
In particular, \citet{ChanLinSaxe2025} combine exact evaluation for a
sampled subset of followers with machine-learning estimates for the
remainder. Our setting differs in the approximation target: rather
than estimating follower or scenario contributions, we predict the
value of graph-structured candidate network configurations.

Another closely related paper is
\citet{SpieckermannMinnerSchiffer2025}, who propose a
reduce-then-optimize framework for the fixed-charge transportation
problem. They use a graph neural network to predict a relevant subset
of variables and then solve a reduced optimization model.
Their approach differs from ours, but it demonstrates how graph-based
learning can be embedded within an exact optimization pipeline for a
network design problem. More broadly, machine learning has been used
to accelerate expensive
optimization components, including column generation
\citep{KraulSeizingerBrunner2023}, while surrogate-assisted
optimization uses regression, interpolation, and ensemble models to
approximate expensive objectives or lower-level value functions
\citep{jin_surrogate-assisted_2011,sinha2016optimistic,sinha2018kriging,dumouchelle2024neurbilo}.
Our strategic problem has this general structure: network topology
defines the candidate solution, while an expensive lower-level
optimizer determines its economic and service outcomes.

Graph-based learning is particularly relevant to network design
because candidate solutions are naturally relational objects: their
value depends not only on the attributes of individual sites but also
on how those sites are connected. Graph neural networks provide a
representation in which node attributes and network structure can be
learned jointly \citep{Battaglia2018,Cappart2021}. The relevant
question in learning-augmented optimization, however, is not simply
whether a graph representation improves predictive accuracy, but
whether it leads the optimization procedure toward better solutions.
This distinction motivates evaluating predictive models through the
exact quality of the network designs they produce when embedded in a
common search procedure.

A related issue arises when surrogate predictions are used to guide
search. Prediction errors that appear small on average can have large
consequences because an optimizer preferentially selects candidates
with attractive predicted values. The search may therefore exploit
optimistic surrogate errors, creating a gap between predictive
accuracy and downstream decision quality. Surrogate-assisted
optimization addresses related concerns by accounting for prediction
uncertainty and by selectively evaluating promising candidates with
the exact model \citep{jin_surrogate-assisted_2011}. Less attention
has been given to
how the scoring rule used within combinatorial search affects
downstream solution quality when several surrogate estimates are available.

Taken together, these literature streams point to two complementary
gaps. Existing spare parts models provide rich representations of
multi-echelon operations but largely treat network structure as
given. At the same time, learning-augmented optimization provides
methods for settings in which exact evaluation is computationally
expensive, but has given limited attention to how surrogate
representation and scoring affect downstream decisions when
predictions guide combinatorial search. We address these gaps through
a framework that combines graph-based surrogate modeling,
conservative scoring, combinatorial search, and exact evaluation of
the final network.

\section{Problem Setting and Decision Structure}
\label{sec:problem-settings}

We consider the design of a two-echelon spare parts inventory system
over a large network of sites. In the incumbent single-echelon
configuration, each site replenishes directly from external suppliers
and manages its inventory independently. In the proposed two-echelon
configuration, a subset of sites is designated as central sites that
replenish inventories for assigned leaf sites, so system performance
depends on both the partition of sites into clusters and the choice
of the central site within each cluster. This induces a hierarchical
decision problem with three interdependent layers: a strategic layer
that determines the network topology, a tactical layer that sets
inventory control parameters for each spare part at each stocking
location, and an operational layer that governs shipment decisions
over time in response to realized demand.

These three decision layers are tightly coupled. Strategic clustering
decisions determine the network topology and therefore the lead
times, pooling opportunities, and admissible material flows across
sites. Those structural choices shape the tactical problem of setting
inventory control parameters for each spare part at each stocking
location. The resulting inventory policies, together with realized
demand, then induce the operational problem of executing material
flows within the network, in particular the timing and quantities of
shipments from central sites to leaf sites. In principle, the
strategic, tactical, and operational layers could be optimized
jointly. In our setting, however, joint optimization is
computationally prohibitive because tactical and operational
decisions must be repeatedly evaluated for each candidate network topology.

The focus of this paper is therefore the strategic design problem:
how to partition the network into clusters and select the central
site of each cluster when the quality of any candidate design can
only be assessed through a computationally expensive lower-level
optimization model. A feasible strategic design consists of a
collection of pairwise-disjoint clusters, each with a designated
central site and a set of leaf sites. Central-to-leaf assignments
must satisfy the operational feasibility conditions imposed by the
business context, most notably distance requirements and a minimum
cluster size. A site assigned to no cluster remains in the incumbent
single-echelon configuration and continues to replenish directly from
external suppliers. Migration to the two-echelon structure is
therefore selective rather than mandatory: the design problem
determines not only how sites should be grouped, but also which sites
are worth grouping at all. The objective is to identify a topology
that improves overall system performance relative to the incumbent
single-echelon configuration while preserving the high service levels
required for maintenance operations.

Because the value of a topology depends on the tactical inventory
policies and operational shipment decisions it induces, the strategic
design problem cannot be reduced to a standard clustering exercise or
assessed from geography alone. Instead, it is a large-scale
combinatorial topology-design problem with expensive lower-level
evaluation: each feasible topology defines an input to the
lower-level model, which returns the economic and service outcomes
used to assess its performance. The next section describes this
lower-level evaluation model and the computational bottleneck it
creates for the learning-augmented solution framework.

\section{Lower-Level Evaluation Model}
\label{sec:nodal-model}

For any fixed network topology, we evaluate its performance using the
two-echelon inventory optimization model developed by
\citet{Caserta2026}. In that model, the network topology is taken as
given, and the optimization determines the tactical inventory
policies and induced operational shipment decisions within the
resulting two-echelon structure under stochastic demand and explicit
business constraints. The role of the lower-level model in the
present paper is therefore purely evaluative: for each candidate
strategic design, it maps the chosen partition of sites and
central-site assignments into exact operational and economic
outcomes. Thus, while \citet{Caserta2026} establishes how to optimize
a given two-echelon configuration, the present paper addresses the
strategic problem of selecting the configuration itself when each
exact evaluation is computationally expensive.

The inputs to the lower-level model include the site assignments
within each cluster, the identity of the central site, part
assortments and demand characteristics, topology-dependent lead
times, and the relevant inventory, stockout, ordering,
transportation, and shipment costs. For a fixed cluster, the model
captures the key distinction between the incumbent single-echelon
system and the proposed two-echelon design: the central site
replenishes from the external supplier, while leaf sites are
replenished from the central site. This changes the effective lead
times and, consequently, the inventory, shipment, and service
outcomes induced by the topology.

We now state the structure of this model explicitly. Consider a
cluster with central site $i$, leaf sites $j \in \SetJ$, spare parts
$k \in \SetK$, periods $t \in \SetT$, and demand scenarios $r \in
\SetS$ obtained from a stochastic demand forecasting model. The
tactical decision variables are the policy parameters: a min level
$\bar{m}^k_l$ for every part $k$ at every site $l$ (central and
leaves) and a max level $\bar{M}^k_i$ for every part at the central
site only, since replenishment at the leaves is governed by the
shipment plan rather than by a local order-up-to level. The
operational decision variables are the order quantities $q^k_{irt}$
placed by the central site with the external supplier, the shipment
quantities $x^k_{ijt}$ from the central site to each leaf, the
resulting inventory and backlog trajectories $s^k_{lrt}$ and
$b^k_{lrt}$, and binary indicators $y^k_{irt}$ and $z^k_{ijt}$ for
order and shipment events, respectively.
Writing $u = (\bar{m},\bar{M})$ for the tactical variables and $v =
(q,x,s,b, y, z)$ for the operational ones, the lower-level problem for a
fixed topology $G$ is the two-stage stochastic program
\begin{equation}
  \label{eq:lower-level}
  Q(G) \;=\; \min_{u,v} \;
  \frac{1}{|\SetS|} \sum_{r \in \SetS}
  \Bigl[
    \underbrace{\textstyle\sum_{k,l,t} c^k_l\, b^k_{lrt}}_{\text{stockout}}
    + \underbrace{\textstyle\sum_{k,l,t} h^k_l\, s^k_{lrt}}_{\text{holding}}
    + \underbrace{\textstyle\sum_{k,t} f^k_i\, y^k_{irt}}_{\text{ordering}}
  \Bigr]
  \;+\;
  \underbrace{\textstyle\sum_{k,j,t} g^k_{ij}\,
  z^k_{ijt}}_{\text{shipment activation}}
\end{equation}
subject to five families of constraints: (i) per-scenario inventory
balance at the central site, which faces its own demand plus the
aggregate dispatch to the leaves and receives supplier orders after
lead time $l^k_i$, and at each leaf, which receives shipments from
the central site after the shorter lead time $l^k_{ij}$; (ii) min/max
policy enforcement at the central site, expressed through
linearized logical constraints stipulating that an order is placed if
and only if the end-of-period inventory position falls to or below
$\bar{m}^k_i$ and, when placed, raises the position to
$\bar{M}^k_i$, with at most one outstanding order per lead-time
window; (iii) leaf-level min maintenance, requiring $s^k_{jrt} \geq
\bar{m}^k_j + 1$ in every scenario and every period after the initial
lead-time window, so that proactive shipments must arrive before any
leaf would breach its min level; (iv) a minimum lot-size constraint
$\bar{M}^k_i - \bar{m}^k_i \geq q^k_{\alpha}$, where $q^k_{\alpha}$
is a prespecified quantile of the empirical distribution of demand
within the supplier lead time; and (v) shipment-consistency
constraints linking quantities to shipment-event indicators, together
with a business rule capping the total number of outgoing shipment
events over the horizon. Uncertainty is handled by Sample Average
Approximation: the scenario set $\SetS$ enters the objective through
the empirical average in~\eqref{eq:lower-level}, and all policy and
shipment constraints must hold for every scenario. The complete
mixed-integer formulation is given in \citet{Caserta2026}.

Service levels are enforced through three complementary mechanisms.
Backlog is penalized in the objective through the stockout cost
$c^k_l$, leaf minimum levels are bounded below by values derived from
a prespecified cycle service level, and the leaf min-maintenance
constraints impose a hard lower bound on inventory after the initial
lead-time window. The service level reported in the computational
study is therefore an output of the exact solution---the fraction of
weeks in which demand is met from available stock, averaged across
scenarios---rather than an explicit model constraint and may
therefore differ slightly from 100\%.

In the implementation used in this paper, the production evaluator
solves formulation~\eqref{eq:lower-level} directly as a multi-item,
multi-period mixed-integer program using FICO Xpress. Although
decomposition approaches have been developed for this problem class
\citep{Caserta2026}, we solve the full formulation because the
evaluator must return exact outcomes for each candidate cluster.
Runtime is driven primarily by the number of demand scenarios, spare
parts, stocking locations, lead-time asymmetries, binary ordering and
shipment decisions, and cluster size. Exact evaluation is feasible
for individual clusters but becomes prohibitively expensive when
embedded in a search over a large number of alternative topologies.

For each candidate design evaluated in our framework, the lower-level
model returns the exact performance measures used at the strategic
layer. In particular, from the exact tactical and operational
solution induced by a topology, we compute the business metrics used
to assess its performance, namely purchase-order spend, average stock value,
and service level, together with their corresponding benchmarks under
the incumbent single-echelon configuration. These exact evaluations
are used in three ways throughout the paper: to label candidate
clusters during data generation, to assess the quality of network
solutions proposed by the surrogate-assisted search, and to
re-evaluate the final selected designs without surrogate
approximation. The next section describes how these exact but
expensive lower-level evaluations are combined with learning and
search to solve the strategic network design problem at scale.

\section{Learning-Augmented Solution Framework}
\label{sec:solution-framework}

This section presents the solution framework for the strategic
network design problem introduced in
Section~\ref{sec:problem-settings}. The key computational challenge
is that evaluating a candidate topology requires solving the exact
lower-level inventory optimization model described in
Section~\ref{sec:nodal-model}. Because this evaluation is expensive,
direct search over the full space of feasible topologies is
impractical. We address this challenge by combining exact evaluation
with graph-based surrogate modeling and metaheuristic search.

The framework separates four roles. First, exact lower-level
optimization provides the labels needed to learn the value of
candidate clusters and the benchmark used to assess all subsequent
decisions. Second, candidate clusters are represented as attributed
graphs and evaluated by a GNN surrogate, allowing the search to
generalize beyond the finite library of exactly evaluated clusters.
Third, an ensemble of GNNs is used to score candidate configurations
conservatively during the search, reducing the influence of
optimistic predictions. Finally, exact set partitioning and
lower-level optimization are used to refine and validate the final
network. The resulting framework uses learning to expand the space
that can be explored while retaining exact optimization for network
selection and final evaluation.

\begin{figure}[t!]
  \centering
  \caption{Structure of the proposed approach decomposed in multiple main
    building blocks. Each block describes a separate step of the
    pipeline that is performed sequentially given the output of the
  previous ones.}
  \label{fig:approach-structure}
  \includegraphics[width=1\textwidth]{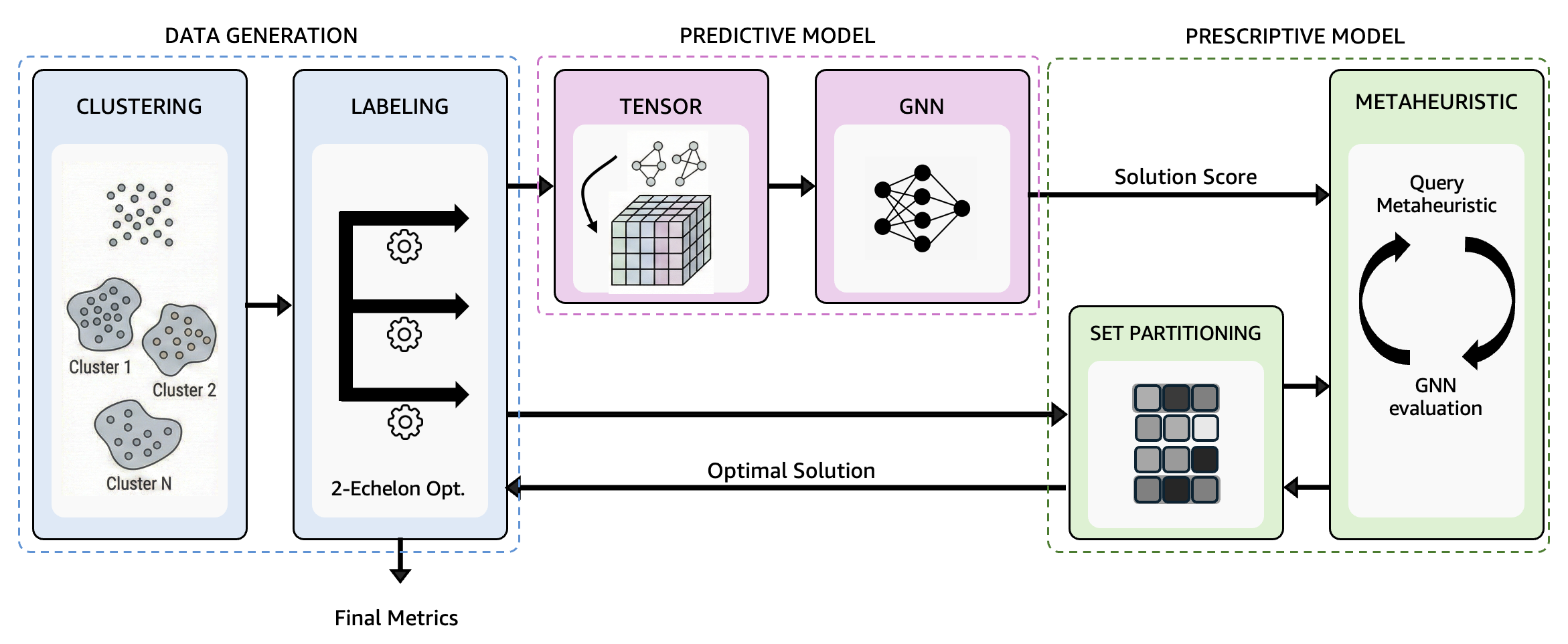}
\end{figure}

Figure~\ref{fig:approach-structure} provides a schematic overview.

\subsection{Exact generation and labeling of candidate clusters}
\label{subsec:cluster-generation-labeling}

We begin by constructing a library of feasible candidate clusters.
Each cluster consists of one designated central site and a set of
leaf sites, and must satisfy the business feasibility requirements of
the application, most importantly the minimum cluster size and the
maximum admissible central--leaf distance. Candidate clusters are
generated by repeatedly solving a clustering problem under
alternative measures of proximity and cluster size configurations,
inducing structural diversity while preserving business feasibility.
Details of the candidate generation procedure are provided in
Section~\ref{ec:candidate-generation} of the appendix.

Each candidate cluster is then evaluated exactly using the
lower-level model of Section~\ref{sec:nodal-model}. For every cluster
$c$, this evaluation returns the exact business metrics under both
the incumbent single-echelon configuration and the proposed
two-echelon configuration. In particular, we denote by
\[
  po_c^t,\quad so_c^t,\quad sl_c^t
\]
the purchase-order spend, average stock value, and service level in
the traditional single-echelon configuration, and by
\[
  po_c^n,\quad so_c^n,\quad sl_c^n
\]
the corresponding quantities under the two-echelon configuration
induced by the cluster, where the superscripts $t$ and $n$ denote the
traditional single-echelon and the new two-echelon configuration,
respectively. These exact evaluations serve two roles in
the framework. First, they define the objective coefficients used to
construct an initial network solution. Second, they provide the
labeled dataset used to train the surrogate model.

Each cluster is therefore labeled using the exact operational model
that defines the strategic objective. The surrogate is not trained on
observational performance data or on a separate approximation of the
operational system, but on outcomes obtained from the lower-level
optimization model.

\subsection{Initial feasible network via set partitioning}
\label{subsec:set-partitioning}

The cluster library generated in the first stage contains only
feasible clusters, but it does not by itself define a full network
design. To assemble a feasible initial network, we solve a set
partitioning model that selects a subset of the evaluated clusters
while avoiding overlap across sites.

Let $\mathcal{I}$ denote the set of sites and let $\mathcal{C}$
denote the set of exactly evaluated candidate clusters. For each
cluster $c \in \mathcal{C}$, let
\[
  a_{ic} =
  \begin{cases}
    1, & \text{if site } i \text{ belongs to cluster } c,\\
    0, & \text{otherwise,}
  \end{cases}
\]
and define the net economic value of cluster $c$ as
\[
  \delta_c = po_c^t + so_c^t - po_c^n - so_c^n.
\]
Thus, $\delta_c$ measures the reduction in purchase-order spend plus
average stock value generated by adopting the two-echelon
configuration represented by cluster $c$ rather than the incumbent
single-echelon benchmark.

The initial network is obtained by solving
\begin{equation}
  \label{eq:spp}
  \max_{x \in \{0,1\}^{|\mathcal{C}|}}
  \left\{
    \sum_{c \in \mathcal{C}} \delta_c x_c :
    \sum_{c \in \mathcal{C}} a_{ic} x_c \le 1,\ \forall i \in \mathcal{I}
  \right\},
\end{equation}
where $x_c=1$ if cluster $c$ is selected and $x_c=0$ otherwise. The
objective maximizes the total estimated savings from the selected
clusters, while the constraints ensure that no site is assigned to
more than one chosen cluster.

Three modeling choices in~\eqref{eq:spp} deserve explicit
justification. First, the coverage constraint is an inequality rather
than an equality because full coverage is neither required nor always
desirable: a site left out of every selected cluster simply remains
in the incumbent single-echelon configuration, replenishing directly
from its external supplier as it does today. Sites therefore need not
join a two-echelon cluster, and an uncovered site contributes zero to
the objective---it generates neither savings nor costs relative to
the single-echelon benchmark against which all savings are measured.
An equality constraint would force the model to include sites in
clusters with negative net value, or render the problem infeasible
whenever the finite library lacks a good cluster for some site.
Second, savings are additive across selected clusters because
clusters are operationally independent by construction: each cluster
is evaluated by the lower-level model as a self-contained system in
which the central site orders from the external supplier and serves
only its own leaves, with no lateral flows between clusters. The
network-level value of a design is therefore the sum of its
cluster-level values, and the network-level performance metrics
reported in the computational study are aggregated from cluster-level
exact evaluations in the same way. Third, the coefficient $\delta_c$
is an exact quantity, not an estimate: every cluster in the library
$\mathcal{C}$ has been evaluated by the lower-level model, so the
baseline produced by~\eqref{eq:spp} is an optimization over exact
evaluations---a strong benchmark against which the surrogate-guided
search must compete.

Model~\eqref{eq:spp} plays two important roles in the overall
framework. First, it provides a strong optimization-based baseline
constructed entirely from exact lower-level evaluations. Second, it
yields the initial feasible solution from which the metaheuristic
search starts. Because the cluster library is only a subset of all
feasible clusters, the solution of~\eqref{eq:spp} may leave some
sites uncovered. Before the search begins, a completion procedure
attempts to raise coverage: uncovered sites that are mutually within
the admissible distance are grouped into provisional clusters via an
auxiliary center-selection model, and only clusters satisfying the
size and center-eligibility rules are retained; remaining sites are attached to
the nearest feasible central site where possible; any site that
cannot be feasibly attached stays in the single-echelon
configuration and can be reconsidered by subsequent neighborhood
moves. The completion procedure is given in
Section~\ref{app:initial-solution} of the appendix.

\subsection{Graph representation of candidate clusters}
\label{subsec:graph-representation}

To generalize beyond the finite library of exactly evaluated
clusters, we represent each candidate cluster as an attributed graph.
This representation is natural because the value of a cluster depends
not only on the attributes of individual sites but also on the
relationships among them, such as proximity, overlap in spare parts
assortments, and the choice of the central site.

Formally, for a cluster $c$ we define a graph
\[
  G_c = (V_c,E_c,V_c^{\mathrm{feat}},E_c^{\mathrm{feat}}),
\]
where $V_c$ is the set of sites in the cluster, $E_c$ is the set of
edges connecting them, $V_c^{\mathrm{feat}}$ collects node features,
and $E_c^{\mathrm{feat}}$ collects edge features. Node features
summarize the intrinsic characteristics of each site, such as site
type, demand-related information, geographic coordinates,
central/leaf status, and aggregate assortment information. Edge
features summarize the relationship between pairs of sites, including
measures of commonality in assortments and cost-related interactions.
This graph construction allows the surrogate to exploit the
relational structure of the strategic design problem rather than
treating each cluster as an unstructured vector of aggregated statistics.

\subsection{Graph neural network surrogate}
\label{subsec:gnn-surrogate}

The goal of the surrogate is to approximate the exact lower-level
evaluation of a candidate cluster at a negligible computational cost.
For a cluster graph $G$, the surrogate maps the graph into predicted
business metrics of the form
\[
  \hat{y}(G) \in \mathbb{R}^{d_{\text{out}}},
\]
where the output components correspond to the lower-level performance
quantities of interest. In our implementation, the surrogate focuses
on the cost-related quantities that drive the strategic objective.

We use a graph-to-vector regressor based on stacked message-passing
layers. Let $H^{(0)} = V^{\mathrm{feat}}$ denote the initial matrix
of node features.
For $\ell=0,\dots,L-1$, the latent node representations are updated according to
\begin{equation}
  \label{eq:gnn-update}
  H^{(\ell+1)}
  =
  \mathrm{SiLU}\!\Big(
    \mathrm{GATv2}^{(\ell)}\!\big(H^{(\ell)},E^{\mathrm{feat}}\big)
  \Big).
\end{equation}
where $\mathrm{GATv2}^{(\ell)}(\cdot)$ denotes an attention-based
message-passing layer with edge conditioning. After $L$ layers, node
embeddings are aggregated into a graph-level representation through
permutation-invariant pooling,
\begin{equation}
  \label{eq:graph-pool}
  h_G = \sum_{i \in V} H_i^{(L)},
\end{equation}
and the final prediction is produced by a multilayer perceptron,
\begin{equation}
  \label{eq:graph-readout}
  \hat{y}(G) = \phi(h_G).
\end{equation}

The graph-based architecture is designed to preserve the structural
information that determines cluster value. Message passing allows the
surrogate to combine site-level characteristics with relationships
among sites, while permutation-invariant pooling ensures that the
prediction does not depend on the ordering of sites in the cluster.
Thus, the model represents a candidate topology as a relational
object rather than as an unordered collection of aggregate
statistics. This representation is particularly important here
because the value of a cluster depends on both the characteristics of
its sites and the way those sites are connected through the proposed
central--leaf structure. Detailed architectural and training choices
are reported in Sections~\ref{app:gnn-architecture}
and~\ref{app:gnn-training} of the appendix.

\subsection{Conservative GNN ensemble}
\label{subsec:pessimistic-ensemble}

A key challenge arises when the surrogate is used inside the search
rather than only for prediction. The search preferentially selects
candidates with high predicted value, so optimistic prediction errors
are more likely to influence subsequent search decisions than
pessimistic errors. A model can therefore have good average
predictive accuracy while still leading the search toward poorly
evaluated configurations.

To obtain a conservative scalar score, we first convert each ensemble
member's four output predictions into a predicted net economic value.
For member $k$, let
$\widehat{po}^{t,(k)}(G)$, $\widehat{so}^{t,(k)}(G)$,
$\widehat{po}^{n,(k)}(G)$, and $\widehat{so}^{n,(k)}(G)$ denote the
predicted traditional and two-echelon PO and stock values, and define
\[
  \widehat{\delta}^{(k)}(G)
  = \widehat{po}^{t,(k)}(G)+\widehat{so}^{t,(k)}(G)
  -\widehat{po}^{n,(k)}(G)-\widehat{so}^{n,(k)}(G).
\]
Instead of using the ensemble mean of these member-level economic
values, the conservative arm scores the candidate using
\begin{equation}
  \label{eq:pessimistic-quantile}
  \widehat{\delta}_{\tau}(G) =
  \mathrm{Quantile}_{\tau}
  \left(
    \widehat{\delta}^{(1)}(G),\dots,\widehat{\delta}^{(K)}(G)
  \right),
\end{equation}
with $\tau<0.5$. Because the search maximizes predicted economic
value, the lower quantile makes a candidate less attractive when
ensemble members do not consistently predict high savings and
therefore reduces the influence of optimistic predictions during
search.

The rationale is prescriptive rather than purely predictive. The
purpose of the ensemble is not only to quantify uncertainty, but to
provide a scoring rule that is better aligned with the way the
optimization algorithm uses predictions. The computational study
therefore compares alternative ensemble scoring rules by the exact
quality of the network designs they produce.

\subsection{Surrogate-assisted search}
\label{subsec:surrogate-assisted-search}

Starting from the initial network produced by the set partitioning
model, we use a variable neighborhood search (VNS) metaheuristic to
explore improved feasible topologies. VNS alternates between
diversification and local improvement, allowing the search to move
across neighborhoods of increasing size
\citep{mladenovic1997variable,hansen2010variable}.

Let $S$ denote a feasible network design and $f(S)$ its
surrogate-based objective value. Starting from an incumbent solution
$S^\star$, the algorithm generates a perturbed solution within the
current neighborhood, applies local improvement, and accepts the
resulting solution when its surrogate score exceeds that of the
incumbent. The search then restarts from the first neighborhood after
an improvement and otherwise moves to the next neighborhood.

Neighborhoods operate directly on cluster topology and include
changes to the central site, transfers of sites across clusters,
cluster splitting, and cluster merging or dissolution when feasible.
These moves preserve the main business constraints while allowing
substantial changes in network structure. The neighborhood structures
used in the computational study are summarized in
Section~\ref{app:vns-neighborhood-summary} of the appendix.

The surrogate makes repeated neighborhood evaluation computationally
feasible. Without it, each candidate move would require an exact
lower-level optimization, whereas the GNN ensemble provides a fast
approximation that allows the search to explore a substantially
larger set of candidate topologies.

\subsection{Final refinement and exact validation}
\label{subsec:exact-validation}

The output of the metaheuristic search is not taken at face value.
Instead, after the search terminates, we aggregate the set of
promising network structures visited during the search and solve one
final set-partitioning recombination problem over this enriched pool. This step
plays the same role as the initial set partitioning model, but now
over a much richer set of candidate structures generated by the
search rather than only the initial library of clusters. For each
experimental arm, the final set partitioning model uses the
corresponding arm-specific surrogate scores as cluster-level objective
coefficients: ensemble means for the mean-scoring arms and
lower-quantile predictions for the conservative-scoring arms. The
model is solved to optimality over the candidate pool generated by
that arm.

Finally, the selected network is re-evaluated exactly using the
lower-level optimization model of Section~\ref{sec:nodal-model}. This
last step is essential because it ensures that the reported
performance of the final solution does not depend on surrogate
approximation. In other words, the surrogate is therefore used to
accelerate exploration, while the
final recommendation is validated using the exact operational model.
The recombination and validation protocol is detailed in
Section~\ref{app:final-sp} of the appendix.

\section{Amazon Case Study: Computational Results}
\label{sec:computational}

\subsection{Experimental setting and evaluation protocol}
\label{subsec:exp_setting}

We evaluate the proposed framework on a large-scale instance derived
from Amazon's North America fulfillment-center network. The case
study contains 246 sites and represents the strategic problem of
designing a two-echelon topology by partitioning sites into feasible
clusters and selecting the corresponding central sites.
Although the present computational study focuses on fulfillment
centers, the framework can also be applied to network settings that
include middle-mile and last-mile facilities.


To construct the dataset for the surrogate model, we generate 3,513
unique feasible candidate clusters using the procedure described in
Section~\ref{subsec:cluster-generation-labeling} and detailed in
Section~\ref{ec:candidate-generation}, and evaluate each cluster
with the exact lower-level optimization model.
Each two-echelon evaluation requires between 10 minutes and 24 hours,
depending on the number of demand scenarios and spare parts in the
cluster. Generating the 3,513 exact cluster evaluations used to train
and evaluate the surrogate required seven days of computation on 10
parallel machines, making repeated exact evaluation impractical
within a large-scale search.
For each cluster, we record the performance of both the
traditional and two-echelon configurations, including purchase-order
(PO) spend, average stock value, and service level.
For each spare part-site pair over a one-year horizon, PO spend is
computed by triggering replenishment whenever inventory reaches the
minimum level and multiplying the order quantity (max level minus
current stock) by the unit price, aggregated across simulated demand
scenarios. Average stock value is approximated as half the difference
between the maximum and minimum inventory levels, multiplied by the
unit price, reflecting the capital tied up in inventory. Service
level is computed across demand scenarios as the fraction of weeks
where the demand can be satisfied using the available stock.
These exact evaluations serve two purposes. First, they provide the
labeled examples used to train the graph-based surrogate. Second,
they define the initial cluster pool used by the constructive set
partitioning model that generates the baseline feasible network for
the subsequent search procedure.

Each evaluated cluster is encoded as a graph whose nodes correspond
to sites and whose edges capture pairwise relationships relevant to
the downstream inventory problem. Node and edge attributes are the
aggregate features introduced in
Section~\ref{sec:solution-framework}. The resulting graph dataset is
split at the cluster level into training, validation, and test sets,
with approximately 80\% of the clusters used for training, 10\% for
validation, and 10\% for testing. The split is performed so that
clusters appearing in the validation and test sets are not used
during training, thereby preventing leakage across identical graph
instances. We train an ensemble of five graph neural networks and use
the validation set for early stopping and model selection.
The same underlying cluster set is used when comparing the GNN with
alternative predictive models, so differences in downstream
performance cannot be attributed to different training data.

The computational study distinguishes clearly between \emph{exact
evaluation} and \emph{surrogate-guided search}. Exact evaluation
refers to solving the lower-level optimization model for a given
cluster or final network design. Surrogate-guided search refers to
the use of the conservative GNN ensemble inside the variable
neighborhood search procedure to score candidate network
modifications rapidly. This distinction is central to the
experimental protocol. The surrogate is used to guide large-scale
exploration of the design space, whereas the performance of the
resulting network designs is assessed using the exact lower-level
optimization model. Predictive performance is evaluated at the
cluster level on held-out observations, while decision quality is
evaluated at the network level through exact re-evaluation of the
solutions produced by the search.

We evaluate the proposed framework against two reference solutions.
The first is the constructive Net-baseline, obtained by solving the
initial set partitioning model over the pool of exactly evaluated
clusters. This provides an optimization-based benchmark constructed
entirely from exact cluster evaluations. The second is a geographic
benchmark, designed to assess how much value can be captured by
clustering sites based on geographic proximity alone. The
Net-baseline serves as the primary reference for measuring
improvements, while the geographic benchmark provides a complementary
comparison with distance-based network design.

\paragraph{Geographic benchmark.} We construct the geographic benchmark by
solving a constrained p-median model with Haversine distance as the
sole clustering criterion. We set p=45 to match the number of
clusters in the constructive Net-baseline. We use the same
eligible-site universe and exclusion rules as in the network-design
procedure and impose the same minimum cluster size m=3,
maximum cluster size of 17, and maximum admissible central–leaf
distance of 500 km. The model jointly selects 45 central sites and
assigns every eligible site to exactly one open central site subject
to these feasibility constraints, minimizing total central–leaf
distance. We solve the resulting mixed-integer program to optimality.
By construction, the geographic network contains 45 clusters and
covers all eligible sites. We then evaluate the resulting topology
using the same exact two-echelon inventory optimization model applied
to all other network designs.

The main outcome measures are PO savings, stock savings, total
savings, service level, and structural characteristics of the
resulting network, such as the number and size distribution of
clusters. Savings are computed relative to the traditional
single-echelon benchmark, so higher values indicate better
performance. Predictive performance is assessed at the cluster level
using held-out observations, while decision quality is assessed at
the network level using exact lower-level evaluation. The
computational study addresses four empirical questions: (i) whether
the GNN provides more accurate predictions than alternative surrogate
representations; (ii) whether conservative ensemble scoring improves
the quality of the decisions produced by the search; (iii) whether
surrogate-guided search improves over the exact constructive baseline
and retains its gains after set-partitioning recombination and exact
validation; and (iv) how well the surrogate generalizes beyond the
distribution represented in the labeled training data and to new
configurations generated during search.

All experiments are run on a cloud machine with 16 CPUs and 64GB of
RAM. The GNN models are implemented in PyTorch, while the lower-level
optimization problems are solved with FICO Xpress. Additional
implementation details for the surrogate and the search are provided
in Sections~\ref{app:gnn-details} and~\ref{app:search-overview} of the
appendix.

\subsection{Surrogate accuracy and conservative scoring}
\label{subsec:surrogate_accuracy}

\begin{figure}[t]
  \centering
  \caption{Comparison of $R^2$ performance across models for PO spend and
  stock level prediction.}
  \label{fig:r2-comparison}
  \includegraphics[width=\linewidth]{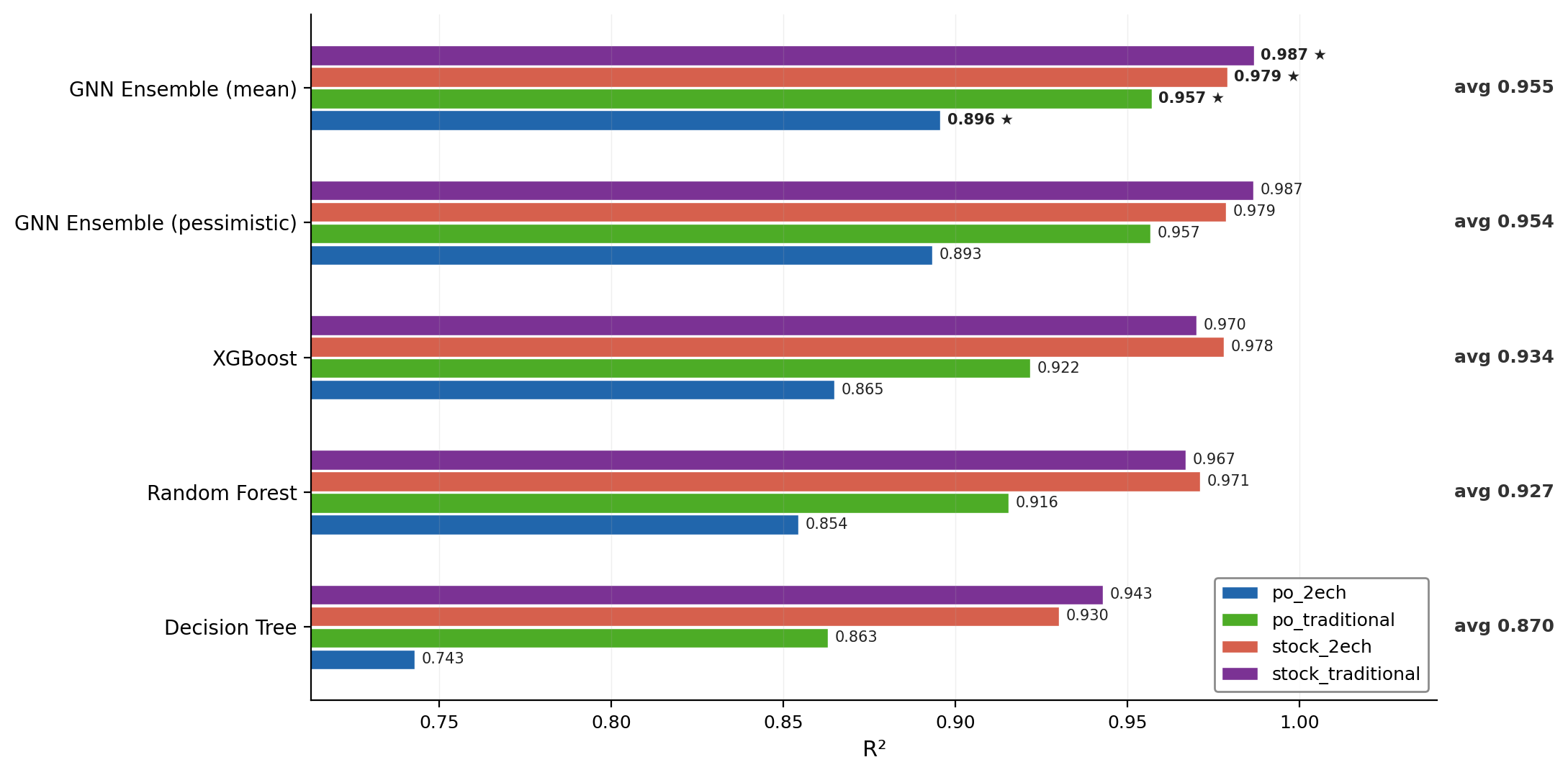}
\end{figure}

We first evaluate the predictive performance of the surrogate on a
held-out test set containing approximately 10\% of the 3{,}513 exactly evaluated
clusters. Figure~\ref{fig:r2-comparison} reports the $R^2$ scores
for the target metrics. Using mean aggregation across the five
ensemble members, the GNN achieves an average $R^2$ of 0.955. The
strongest non-graph alternative, XGBoost, achieves 0.934, followed
by Random Forest at 0.927 and a single Decision Tree at 0.870. The
GNN therefore provides the most accurate predictions across the
models considered, improving average $R^2$ by approximately two
percentage points over the strongest tabular alternative.

These results establish predictive performance on clusters drawn from
the same underlying distribution as the training data, but they do
not establish that the surrogate will remain accurate on
configurations encountered during optimization. The search
systematically explores new cluster compositions and may therefore
move away from the distribution represented in the labeled library.
Section~\ref{subsec:distribution_shift} examines this issue through
controlled distribution shifts.

\begin{table}[ht]
  \caption{Predictive models architectures for cluster cost estimation.}
  \label{tab:predictive-models}
  \centering
  \resizebox{\textwidth}{!}{%
    \begin{tabular}{llccp{5.5cm}p{4.5cm}}
      \toprule
      \textbf{Model} & \textbf{Type} & \textbf{Input
      Representation} & \textbf{Nb.\ Features} & \textbf{Key
      Hyperparameters}\\
      \midrule
      Decision Tree
      & Single regression tree
      & Tabular
      & 91
      & max\_depth=10\\
      \midrule
      Random Forest
      & Bagged tree ensemble
      & Tabular
      & 91
      & 100 estimators, max\_depth=15\\
      \midrule
      XGBoost
      & Boosted tree ensemble
      & Tabular
      & 91
      & 200 estimators, max\_depth=8, lr=0.1\\
      \midrule
      GNN Ensemble
      & Ensemble of 5 GATv2
      & Graph Tensor
      & 11 node + 3 edge
      & 4 GATv2 layers (2/2/1/1 heads), hidden=32, dropout=0.2,
      SiLU, AdamW + OneCycleLR, patience=50\\
      \bottomrule
  \end{tabular}}
\end{table}

\begin{figure}[t]
  \centering
  \caption{Distribution of predictions from the GNN ensemble across the five
    model instances. Each box-plot shows the distribution of
    predictions for a single test sample. Samples are sorted by actual
    value along the x-axis, while the y-axis shows the prediction error
  in percentage.}
  \label{fig:gnn_distribution}
  \includegraphics[width=0.8\linewidth]{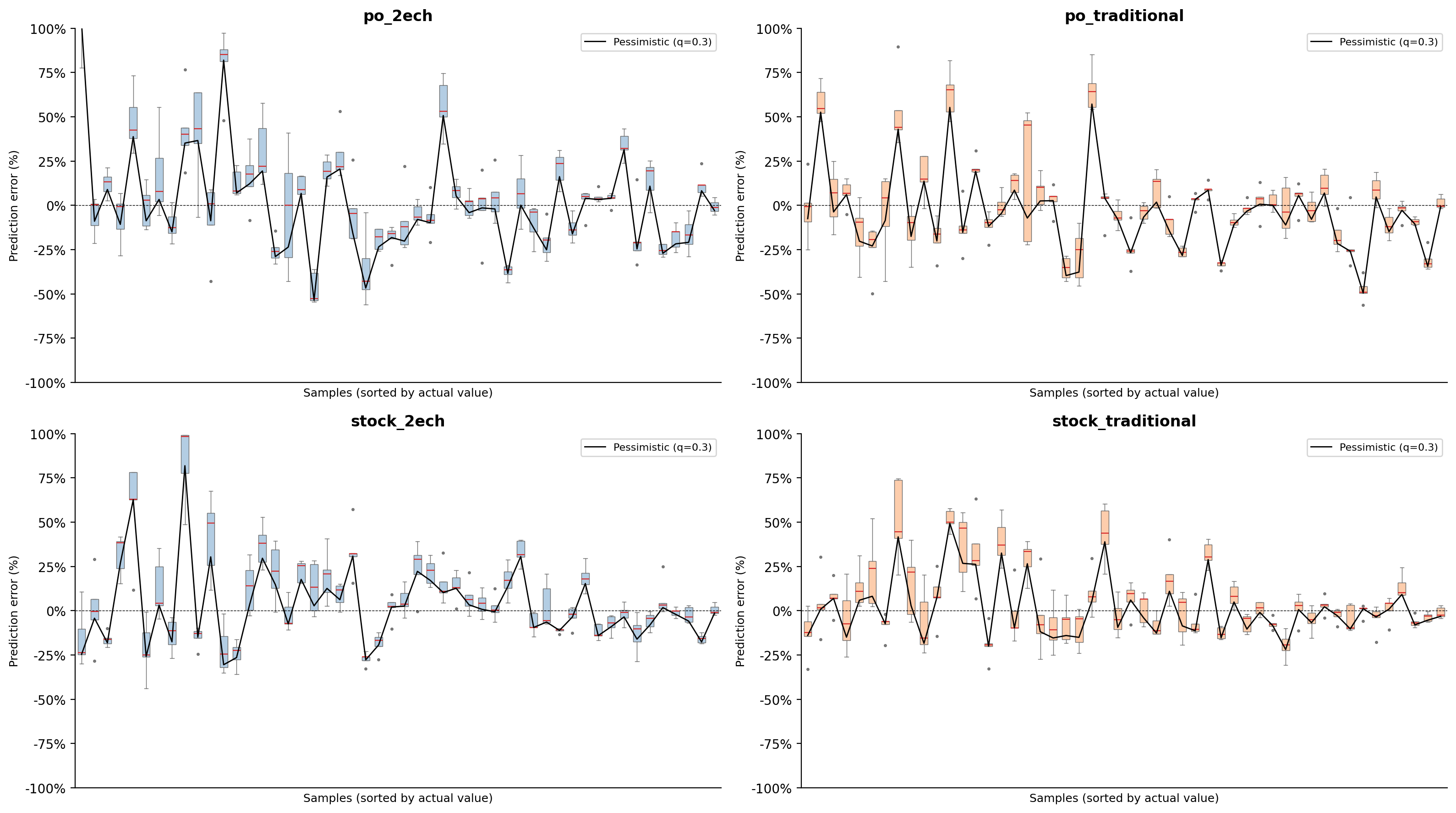}
\end{figure}

We next examine the effect of conservative ensemble scoring. The
search uses the 30th percentile of the five GNN predictions rather
than their mean, as described in
Section~\ref{subsec:pessimistic-ensemble}. On the held-out test set,
this change has almost no effect on predictive accuracy: average
$R^2$ decreases only from 0.955 under mean aggregation to 0.954 under
conservative aggregation. Figure~\ref{fig:gnn_distribution}
illustrates how the ensemble predictions vary across test clusters
and how the lower-quantile rule shifts the score toward the
conservative side of the prediction distribution. Whether this change
in scoring improves the quality of the network designs produced by
the search is examined directly in Section~\ref{subsec:component_analysis}.

\subsection{Generalization Under Distribution Shift}
\label{subsec:distribution_shift}

The predictive results in Section~\ref{subsec:surrogate_accuracy} are
based on a random split of the labeled cluster library. Such a split
measures interpolation within the observed data distribution, whereas
the surrogate is ultimately used to evaluate new configurations
generated by the search. We therefore examine how predictive
performance changes under several forms of distribution shift.

We retrain the GNN ensemble using grouped train--test splits while
holding the architecture, training procedure, and ensemble size
fixed. We consider three types of holdout. First, a size-based split
withholds all clusters containing at least eight sites, requiring the
surrogate to extrapolate to larger clusters than those observed
during training. Second, a central-site split withholds all clusters
associated with 19 randomly selected central sites. Third, geographic
splits withhold all clusters whose central site belongs to a given
sub-region. The random train-test experiment of
Section~\ref{subsec:surrogate_accuracy} provides the reference case.

\begin{table}[t]
  \caption{Generalization of the surrogate ensemble under grouped
    splits. Each row is an independent retraining of the five-member
    ensemble with the stated clusters withheld; columns report the
    average $R^2$ over the four reported targets under mean
    aggregation and under the conservative $\tau{=}0.3$ aggregation,
    followed by the mean-aggregated $R^2$ of the two savings drivers
  taken separately.}
  \label{tab:distribution_shift}
  \centering
  \begin{tabular}{lrrrrrr}
    \toprule
    & & & \multicolumn{2}{c}{Avg.\ $R^2$} &
    \multicolumn{2}{c}{Mean agg., by target} \\
    \cmidrule(lr){4-5} \cmidrule(lr){6-7}
    \textbf{Withheld group} & \textbf{Development} & \textbf{Test} &
    \textbf{Mean} & \textbf{$\tau{=}0.3$} &
    \textbf{PO} & \textbf{Stock} \\
    \midrule
    None (random $80/10/10$)        & $3{,}161$ & $352$   & $0.955$ &
    $0.954$ & $0.877$ & $0.982$ \\
    \midrule
    Clusters with ${\geq}8$ sites & $2{,}419$ & $1{,}094$ & $0.919$
    & $0.876$ & $0.853$ & $0.952$ \\
    $19$ central sites            & $3{,}062$ & $451$   & $0.802$ &
    $0.816$ & $0.605$ & $0.907$ \\
    \midrule
    Sub-region: South Central     & $2{,}788$ & $725$   & $0.785$ &
    $0.742$ & $0.674$ & $0.820$ \\
    Sub-region: Northeast         & $2{,}489$ & $1{,}024$ & $0.769$
    & $0.735$ & $0.578$ & $0.898$ \\
    Sub-region: West              & $3{,}126$ & $387$   & $0.706$ &
    $0.501$ & $0.589$ & $0.777$ \\
    Sub-region: Central           & $2{,}654$ & $859$   & $0.465$ &
    $0.528$ & $-0.263$ & $0.856$ \\
    Sub-region: Southeast         & $2{,}995$ & $518$   & $0.343$ &
    $0.498$ & $-0.209$ & $0.835$ \\
    \bottomrule
  \end{tabular}
\end{table}

Table~\ref{tab:distribution_shift} reveals a substantial difference
between interpolation and extrapolation. Under the random split, the
ensemble achieves an average $R^2$ of 0.955. Withholding large
clusters reduces average $R^2$ only moderately, to 0.919, despite
requiring predictions for cluster sizes absent from the training set.
Withholding central sites has a larger effect, reducing $R^2$ to
0.802. Geographic holdouts are substantially more challenging, with
average $R^2$ ranging from 0.343 to 0.785 across sub-regions. These
results show that the surrogate generalizes relatively well across
cluster size but is more sensitive to geographic shifts in the
composition of the training data.

The two economic targets also respond differently to geographic
shift. Stock-value predictions remain comparatively stable, with $R^2$
no lower than 0.777 across the geographic holdouts, whereas PO-spend
predictions deteriorate sharply in some regions, reaching negative $R^2$
in the Central and Southeast holdouts. These negative values indicate
that the surrogate predicts PO spend less accurately than a model
that simply predicts the mean PO spend of the corresponding held-out
observations as constant predictions. This difference suggests that
some components of
cluster value are more transferable across regions than others. In
particular, the results indicate that geographic coverage of the
training data is important when the surrogate is used to predict
purchasing-related outcomes.

Overall, the results show that generalization depends strongly on the
type of distribution shift. The GNN remains relatively accurate when
extrapolating to larger clusters, but performance deteriorates more
substantially when central sites or entire geographic regions are
absent from the training data. Geographic shifts are particularly
consequential for PO-spend predictions. These results highlight the
importance of geographic coverage in the training data when the
surrogate is used across heterogeneous regions.

The need for generalization beyond the labeled library is
particularly relevant during search. In the production run using the
conservative $\tau=0.3$ score, the search visits 18,271 distinct cluster
compositions, of which only 48 (0.3\%) appear in the labeled library.
The surrogate therefore evaluates primarily configurations that
have not previously been solved by the lower-level model.
Evaluating all of these configurations exactly would require more
than 3,000 machine-hours even under the optimistic assumption of
10 minutes per evaluation. Surrogate scoring thus makes it
possible to explore a much larger set of candidate configurations
than could practically be evaluated with the exact model.
Section~\ref{subsec:search_generated_validation}
examines how well the surrogate performs on these
search-generated candidates.

\subsection{Surrogate performance on search-generated candidates}
\label{subsec:search_generated_validation}

The grouped holdout experiments in Section~\ref{subsec:distribution_shift}
provide controlled tests of distribution shift, but they do not directly
measure surrogate performance on the configurations generated by the
optimization procedure. This distinction is important because the search
does not sample candidate clusters uniformly: it preferentially explores
configurations with attractive predicted values. We therefore conduct an
additional validation experiment on previously unseen clusters visited by
the variable neighborhood search.

From the set of visited clusters that do not belong to the labeled library,
we draw a stratified sample of 300 configurations. The sample contains 100
randomly selected visited clusters, 100 clusters drawn from the top 20\% of
surrogate scores, and 100 clusters drawn from the top 5\%. Each sampled
cluster is then evaluated using the exact lower-level optimization model.
This design allows us to assess not only average predictive performance
outside the labeled library, but also performance in the high-value tail
that is most consequential for the search.

For each cluster $c$, we define the signed percentage prediction error as
\begin{equation}
  e_c =
  100\frac{\widehat{V}_c-V_c}{|V_c|},
  \label{eq:signed_prediction_error}
\end{equation}
where $\widehat{V}_c$ is the surrogate value and $V_c$ is the corresponding
exact value. Because the search maximizes economic value, positive
errors indicate optimistic predictions, for which the surrogate makes
a candidate appear more attractive than its exact evaluation, whereas
negative errors indicate conservative predictions.

\begin{table}[t]
  \caption{Predictive performance on previously unseen clusters generated
    during search. Search-generated clusters are stratified according to their
    surrogate-predicted value. All performance measures are computed against
  exact lower-level evaluations.}
  \label{tab:search_generated_accuracy}
  \centering
  \begin{tabular}{lccccc}
    \toprule
    \textbf{Search-generated sample}
    & \textbf{$n$}
    & \textbf{$R^2$}
    & \textbf{Spearman $\rho$}
    & \textbf{Mean Signed error (\%)}
    & \textbf{Overest.} \\
    \midrule
    Random visited clusters
    & 100 & 0.88 & 0.93 & $-1.8\%$ & $31\%$ \\
    Top 20\% by surrogate score
    & 100 & 0.81 & 0.89 & $-2.6\%$ & $27\%$ \\
    Top 5\% by surrogate score
    & 100 & 0.74 & 0.85 & $-3.1\%$ & $24\%$ \\
    \midrule
    All sampled clusters
    & 300 & 0.84 & 0.91 & $-2.4\%$ & $27\%$ \\
    \bottomrule
  \end{tabular}
\end{table}

Table~\ref{tab:search_generated_accuracy} reports the results for the
conservative GNN ensemble. Across all 300 search-generated clusters, the
surrogate achieves an $R^2$ of 0.84 and a Spearman rank correlation of
0.91. Predictive accuracy decreases as attention shifts toward the
candidates most favored by the search: $R^2$ falls from 0.88 for randomly
visited clusters to 0.74 among the top 5\% of surrogate scores. Rank
correlation, however, remains comparatively high, declining from 0.93 to
0.85. This distinction is important because the search relies primarily on
relative comparisons among candidate configurations rather than on perfectly
calibrated predictions of their absolute values.

The conservative ensemble also remains cautious in the high-value
tail. Mean signed error is negative in all three strata, and only 24\%
of the top-5\% candidates receive optimistic predictions. Thus,
although absolute predictive accuracy deteriorates for the extreme
configurations preferred by the optimizer, the surrogate continues
to rank these candidates effectively while tending to assign values
below their exact evaluations.

\paragraph{Graph representation under search-induced shift.}
We next examine whether the predictive advantage of the GNN persists
on configurations generated by the search. We compare the GNN with
the strongest tabular benchmark, XGBoost, using the same exactly
evaluated search-generated clusters. Holding the candidate set fixed
allows us to compare how well the two surrogate specifications
predict and rank configurations encountered during search. Because
candidate ranking is central to the search, we report both $R^2$ and
Spearman rank correlation.

\begin{table}[t]
  \caption{Predictive performance of GNN and tabular surrogates on
    search-generated clusters. Both models are evaluated against exact
  lower-level outcomes on the same candidate configurations.}
  \label{tab:search_generated_models}
  \centering
  \begin{tabular}{lcccc}
    \toprule
    &
    \multicolumn{2}{c}{\textbf{GNN}}
    &
    \multicolumn{2}{c}{\textbf{XGBoost}} \\
    \cmidrule(lr){2-3}
    \cmidrule(lr){4-5}
    \textbf{Search-generated sample}
    & \textbf{$R^2$}
    & \textbf{Spearman $\rho$}
    & \textbf{$R^2$}
    & \textbf{Spearman $\rho$} \\
    \midrule
    Random visited clusters
    & 0.88 & 0.93
    & 0.82 & 0.87 \\
    Top 20\% by surrogate score
    & 0.81 & 0.89
    & 0.67 & 0.72 \\
    Top 5\% by surrogate score
    & 0.74 & 0.85
    & 0.48 & 0.55 \\
    \bottomrule
  \end{tabular}
\end{table}

Table~\ref{tab:search_generated_models} shows that the GNN maintains
stronger predictive and ranking
performance than XGBoost across all three samples, with the
difference becoming more pronounced among candidates receiving the
highest surrogate scores. For randomly visited clusters, the GNN
achieves an $R^2$ of 0.88 and a Spearman rank correlation of 0.93,
compared with 0.82 and 0.87 for XGBoost. Among candidates sampled
from the top 5\% of surrogate scores, the GNN retains an $R^2$ of 0.74
and a rank correlation of 0.85, whereas the corresponding values
for XGBoost fall to 0.48 and 0.55. The difference in rank
correlation is particularly relevant because the search relies on
the relative ordering of candidate configurations.

This result provides a different perspective from the conventional
held-out comparison in Section~\ref{subsec:surrogate_accuracy}. The
approximately two-percentage-point advantage of the GNN under a random
train--test split understates the difference between the models in the
region of the design space that is most consequential for optimization.
The graph representation appears particularly valuable for preserving the
relative ordering of attractive configurations as the search moves away
from the labeled cluster library. Together with the equal-budget
component analysis
in Section~\ref{subsec:component_analysis}, this analysis allows us to
distinguish whether the graph representation merely improves conventional
predictive accuracy or also improves the information available to the
optimization procedure.

\paragraph{Optimizer-induced overestimation and conservative scoring.}
Finally, we use the search-generated sample to examine the mechanism
motivating conservative ensemble scoring. For each exactly evaluated
cluster, we compute both the mean ensemble prediction and the conservative
$\tau=0.3$ score. The candidate set is therefore held fixed, and only the
rule used to aggregate the five GNN predictions changes.

\begin{table}[t]
  \caption{Effect of ensemble scoring on prediction error for
    search-generated clusters. Mean signed error (mse) and the conservative scores are evaluated on
  the same clusters using exact lower-level outcomes as the benchmark.}
  \label{tab:search_generated_quantile}
  \centering
  \begin{tabular}{lcccc}
    \toprule
    &
    \multicolumn{2}{c}{\textbf{Mean aggregation}}
    &
    \multicolumn{2}{c}{\textbf{Conservative $\tau=0.3$}} \\
    \cmidrule(lr){2-3}
    \cmidrule(lr){4-5}
    \textbf{Search-generated sample}
    & \textbf{mse (\%)}
    & \textbf{Overest.}
    & \textbf{mse (\%)}
    & \textbf{Overest.} \\
    \midrule
    Random visited clusters
    & $+0.4\%$ & 48\%
    & $-1.8\%$ & 31\% \\
    Top 20\% by surrogate score
    & $+4.7\%$ & 63\%
    & $-2.6\%$ & 27\% \\
    Top 5\% by surrogate score
    & $+9.8\%$ & 78\%
    & $-3.1\%$ & 24\% \\
    \bottomrule
  \end{tabular}
\end{table}

Table~\ref{tab:search_generated_quantile} reveals a clear selection effect
under mean aggregation. Among randomly selected visited clusters, prediction
errors are approximately centered, with a mean signed error of $+0.4\%$ and
an overestimation rate of 48\%. Among candidates receiving higher
conservative GNN scores, however, optimistic errors under mean
aggregation become increasingly
prevalent. The overestimation rate rises to 63\% among the top 20\%
sample and to 78\% among the top 5\% sample, where the mean predictition exceed
exact the value by 9.8\% on average. hus, prediction errors that are
approximately symmetric among randomly visited clusters become
strongly skewed toward overestimation among the high-scoring
candidates considered here.

Conservative aggregation largely removes this pattern. Using the 30th
percentile, mean signed error remains negative across all three
samples, while the overestimation rate decreases as attention shifts
toward candidates receiving higher conservative scores. Among the
top-5\% sample, only 24\% of clusters are overestimated under
conservative aggregation, compared with 78\% when the same clusters
are evaluated using the ensemble mean. Because the candidate set
and ensemble predictions are held fixed, this comparison provides
direct evidence for the mechanism motivating conservative scoring:
candidates that appear attractive under surrogate guidance can
exhibit substantial optimistic error under mean aggregation, while
lower-quantile aggregation reduces the influence of such errors.

Taken together, the search-generated validation provides two insights
that are not visible from conventional held-out prediction metrics
alone. First, the GNN maintains stronger predictive and ranking
performance than the tabular alternative on the common set of
high-scoring search-generated candidates considered here. Second,
holding the GNN ensemble and candidate set fixed, conservative
aggregation substantially reduces optimistic prediction errors among
these candidates. These results provide evidence on how surrogate
specification and scoring behave in the portion of the design space
explored by the search. Section~\ref{subsec:component_analysis} then
examines the more
consequential question of whether these differences translate into
better network designs under exact downstream evaluation.

\subsection{Component analysis: what drives decision quality}
\label{subsec:component_analysis}

To isolate the contribution of the principal methodological choices,
we conduct an equal-budget component analysis on the North America
instance. All experimental arms use the same candidate-cluster
labels, initial set-partitioning baseline,
variable-neighborhood-search budget, set-partitioning recombination
procedure, and final lower-level validation. They differ only in the
surrogate representation and the rule used to aggregate ensemble
predictions. Because the search is stochastic, each principal arm is
run independently six times, and every resulting network is evaluated
using the exact lower-level optimization model. We report arm-level
means and standard deviations across the six independent runs and
interpret differences across arms as descriptive case-study evidence.
All reported improvements are relative to the exact Net-baseline.

For the XGBoost arms, we construct an ensemble of five models trained
on independently resampled versions of the training data, paralleling
the five-member GNN ensemble. This allows us to apply the same
aggregation rules to both surrogate classes and separate, as far as
possible, the effect of surrogate choice from the effect of scoring.
For both ensembles, the mean-scoring arm uses the mean of the five
member-level predictions, whereas the conservative arm uses their
30th percentile.

Table~\ref{tab:component_analysis} show descriptive differences in
mean downstream decision quality across the component arms. Under
mean scoring, the GNN achieves mean exact total savings of 25.4\%
above the Net-baseline, compared with 17.4\% for XGBoost, a
difference of 8.0 percentage points across six independent runs per
arm. Under conservative scoring after set-partitioning
recombination, the corresponding improvements are 30.5\% for the GNN
and 21.9\% for XGBoost, a difference of 8.6 percentage points.
Within this case study, the predictive advantage of the GNN
documented in Section~\ref{subsec:surrogate_accuracy} is therefore
accompanied by higher mean
exact savings when the surrogate is embedded in the search.

Conservative scoring provides a second comparison. Holding the
surrogate class fixed, replacing mean aggregation with the
lower-quantile score increases mean exact total savings by 4.5
percentage points for XGBoost and by 5.1 percentage points for the
GNN. The direction of the difference is therefore the same for both
surrogate classes: in these runs, conservative scoring produces
networks with higher mean exact savings than mean aggregation. These
differences are descriptive, but their direction is consistent across
both surrogate classes.

Table~\ref{tab:component_mechanism} provides
additional evidence on the mechanisms underlying the differences in
downstream performance. Under mean aggregation, a large share of the
clusters selected into the resulting networks are overestimated by
the surrogate: 76\% for XGBoost and 68\% for the GNN. Conservative
scoring reduces these frequencies to 39\% and 28\%, respectively, and
shifts mean signed error from positive to slightly negative. This
pattern is consistent with the search-generated analysis in
Table~\ref{tab:search_generated_quantile}
and indicates that conservative scoring reduces the prevalence of
optimistic predictions among clusters that ultimately enter the
network. At the same time, the GNN exhibits higher rank correlation
with exact cluster values than XGBoost under both scoring rules.
These results suggest complementary roles for the two
methodological choices: the GNN provides more reliable rankings of
selected configurations, while conservative scoring reduces the
influence of optimistic prediction errors.

The GNN with conservative scoring achieves the highest mean exact
savings among the component arms. Taken together, Tables 6 and 7
suggest complementary roles for surrogate choice and scoring: the GNN
is associated with stronger rankings of selected configurations,
while conservative aggregation reduces the prevalence of optimistic
prediction errors. The final set-partitioning recombination provides
an additional 3.1-percentage-point improvement, increasing mean exact
total savings from 27.4\% before recombination to 30.5\% after recombination.
Thus, the results suggest complementary roles for the main elements
of the framework: the GNN provides stronger rankings, conservative
scoring limits the influence of optimistic predictions, and
set partitioning recombination extracts additional value from the
candidate pool generated during search.

\begin{table}[t]
  \caption{Component analysis on the North America instance. PO and stock
    are six-run mean improvements relative to the exact Net-baseline; Total
    is the mean combined-savings improvement $\pm$ standard deviation.
    Principal surrogate arms include final set-partitioning
    recombination using their arm-specific scores; the pre-SP result and
  single-evaluation geographic benchmark are shown separately.}
  \label{tab:component_analysis}
  \centering
  \begin{tabular}{lrrr}
    \toprule
    \textbf{Arm}
    & \textbf{PO}
    & \textbf{Stock}
    & \textbf{Total} \\
    \midrule
    Net-baseline
    & Ref. & Ref. & Ref. \\

    Geography ($p$-median)
    & $-51.4\%$ & $-48.8\%$ & $-49.9\%$ \\
    \midrule

    XGBoost + mean
    & $+8.9\%$ & $+23.7\%$
    & $+17.4\% \pm 2.7~\mathrm{pp}$ \\

    XGBoost + conservative
    & $+12.1\%$ & $+29.8\%$
    & $+21.9\% \pm 2.4~\mathrm{pp}$ \\

    GNN + mean
    & $+14.8\%$ & $+34.1\%$
    & $+25.4\% \pm 2.2~\mathrm{pp}$ \\

    GNN + conservative (post-SP)
    & $+18.7\%$ & $+40.2\%$
    & $+30.5\% \pm 1.9~\mathrm{pp}$ \\
    \midrule

    GNN + conservative (pre-SP)
    & $+16.3\%$ & $+35.7\%$
    & $+27.4\% \pm 2.0~\mathrm{pp}$ \\
    \bottomrule
  \end{tabular}
\end{table}

The geographic benchmark provides a complementary comparison. The
distance-based p-median solution yields an economic objective 49.9\%
below the Net-baseline, indicating that geographic proximity alone
does not necessarily identify clusters with attractive inventory
and purchasing outcomes. This comparison should be interpreted
cautiously, however, because the geographic benchmark requires full
coverage of all eligible sites, whereas the Net-baseline allows
selective participation in the two-echelon network and may
subsequently increase coverage through the completion procedure.

\begin{table}[t]
  \caption{Surrogate behavior on clusters selected by the search. Values are
    computed against exact lower-level evaluations of the clusters included in
  the resulting networks.}
  \label{tab:component_mechanism}
  \centering
  \begin{tabular}{lccc}
    \toprule
    \textbf{Arm}
    & \textbf{Overest.}
    & \textbf{Mean Signed error (\%)}
    & \textbf{Spearman $\rho$} \\
    \midrule
    XGBoost + mean
    & $76\%$ & $+8.7\%$ & 0.61 \\

    XGBoost + conservative
    & $39\%$ & $-1.9\%$ & 0.69 \\

    GNN + mean
    & $68\%$ & $+6.1\%$ & 0.78 \\

    GNN + conservative
    & $28\%$ & $-2.4\%$ & 0.86 \\
    \bottomrule
  \end{tabular}
\end{table}

\subsection{End-to-end network performance: overall results}
\label{subsec:validated_performance}

Table~\ref{tab:results_m1} provides a consolidated view of the
end-to-end performance of the proposed framework. Although some of
the economic results have been reported in the component analysis of
the previous section, we collect them here together with service and
network-structure measures to provide an overall assessment of the
final solution. The table compares the exact constructive baseline,
the network identified by the GNN-guided search with conservative
ensemble scoring, and the final network obtained after
set-partitioning recombination. All reported performance measures are
based on exact evaluation using the two-echelon lower-level model.

The end-to-end results show substantial economic improvements over
the constructive Net-baseline. The GNN-guided search with
conservative scoring improves mean exact total savings by 27.4\%, and
final set-partitioning recombination increases this improvement to
30.5\%. The gains reflect improvements in both purchase-order spend
and stock value, as reported in Table~\ref{tab:results_m1}.

These economic gains are achieved while maintaining an average
service level of approximately 99.8\%. The final recombination also
modifies the network structure, increasing the mean number of
clusters from 27 to 28 and reducing the maximum cluster size from
17 to 14, while preserving the economic improvements identified during search.

\begin{table}[t]
  \caption{End-to-end network performance of the proposed framework on the
    North America instance. Savings are reported as percentage improvement
    relative to the exact Net-baseline. Values are mean $\pm$ standard
    deviation across six independent VNS replications. All networks are
  evaluated using the exact two-echelon lower-level model.}
  \label{tab:results_m1}
  \centering
  \footnotesize
  \setlength{\tabcolsep}{3pt}
  \begin{tabular}{lrrrrrr}
    \toprule
    & \multicolumn{3}{c}{Savings Improvement vs.\ Net-baseline}
    & \multicolumn{3}{c}{Network Structure} \\
    \cmidrule(lr){2-4} \cmidrule(lr){5-7}
    \textbf{Solution}
    & \textbf{PO}
    & \textbf{Stock}
    & \textbf{Total}
    & \textbf{\# Clusters}
    & \textbf{Min}
    & \textbf{Max} \\
    \midrule
    Net-baseline
    & Ref.
    & Ref.
    & Ref.
    & 45
    & 3
    & 10 \\

    GNN + conservative VNS
    & $+16.3\% \pm 1.8~\mathrm{pp}$
    & $+35.7\% \pm 2.1~\mathrm{pp}$
    & $+27.4\% \pm 2.0~\mathrm{pp}$
    & 27
    & 3
    & 17 \\

    GNN + conservative VNS + final SP
    & $+18.7\% \pm 1.7~\mathrm{pp}$
    & $+40.2\% \pm 1.8~\mathrm{pp}$
    & $+30.5\% \pm 1.9~\mathrm{pp}$
    & 28
    & 3
    & 14 \\
    \bottomrule
  \end{tabular}
\end{table}

\begin{figure}[t!]
  \centering
  \caption{Final validated two-echelon network identified by the proposed
    framework for the 246-site North America fulfillment-center instance.
    Stars denote selected central sites, circles denote assigned leaf sites,
    connecting lines represent central-to-leaf assignments, and colors
  identify clusters.}
  \label{fig:map}
  \includegraphics[width=0.8\textwidth]{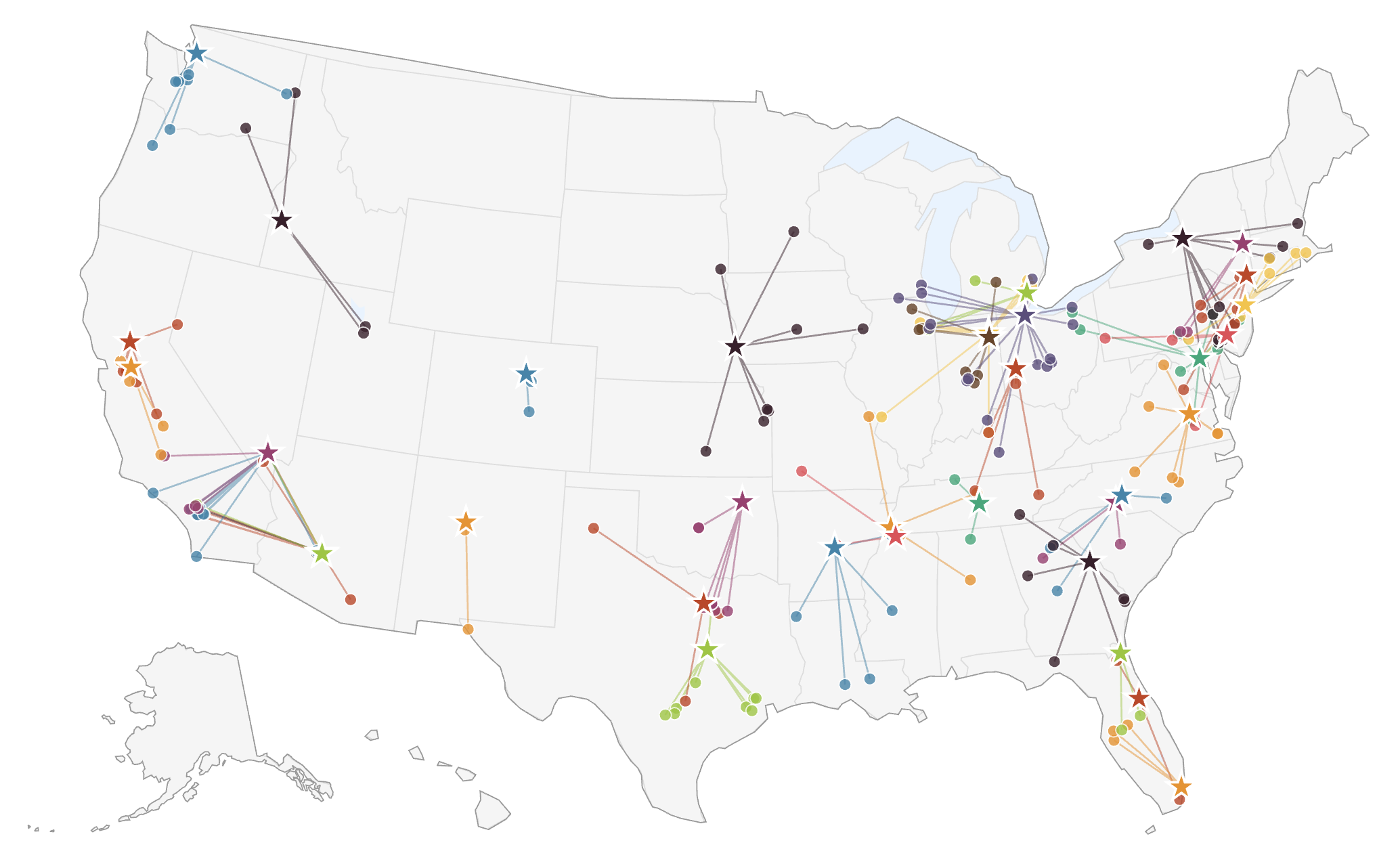}
\end{figure}

Figure~\ref{fig:map} illustrates the final validated network. The
resulting topology is heterogeneous: some regions consolidate around
relatively large central sites, whereas others retain smaller, more
localized clusters. The framework therefore allows the degree of
centralization to vary across regions rather than imposing a uniform
network structure. Additional descriptive evidence on the convergence
and dispersion of the VNS across independent restarts is reported in
Section \ref{app:vns-convergence} of the appendix.

\section{Discussion and Conclusions}
\label{sec:conclusions}

This paper develops a learning-augmented optimization framework for
the strategic design of large-scale two-echelon spare parts networks,
combining exact lower-level evaluation with surrogate-guided search
and final network refinement. In the Amazon case study, the framework
improves combined savings by 30.5\% relative to an optimization-based
baseline constructed from exact cluster evaluations, while
maintaining an average service level of approximately 99.8\%. Beyond
this overall improvement, the component analysis highlights two
broader findings: the choice of surrogate can affect
downstream decision quality, and the way surrogate predictions are
used within optimization matters when the search preferentially
explores candidates with favorable prediction errors.

\subsection{Managerial insights}
\label{subsec:managerial-insights}

Our results provide four main implications for firms redesigning large-scale
spare parts networks.

\textbf{Insight 1: Strategic network design should be driven by operational
value rather than geography alone.}
Geographic proximity is an important feasibility consideration in
two-echelon network design, but it does not capture the inventory and
purchasing consequences of a candidate topology. The value of
grouping a set of sites depends on the inventory policies,
replenishment decisions, lead times, and shipment flows induced by
that structure. In the case study, the distance-based p-median
benchmark performs substantially worse economically than the
Net-baseline. The managerial implication is that geographic
clustering can provide a useful starting point or feasibility
criterion, but strategic network decisions should ultimately account
for the operational and economic consequences of the resulting topology.

\textbf{Insight 2: The appropriate degree of centralization can vary across the
network.}
The final validated network in Figure~\ref{fig:map} is structurally
heterogeneous.
Some regions consolidate around relatively large central sites serving
several downstream locations, whereas others retain smaller, more localized
clusters. The resulting design therefore differs from policies that impose
a common cluster size or degree of centralization throughout the network.
For managers, this suggests that two-echelon redesign should allow the
network structure to adapt to local operating conditions rather than impose
a uniform template across geographically and operationally heterogeneous
regions.

\textbf{Insight 3: Learning components should be evaluated by the decisions
they produce, not by predictive accuracy alone.}
When machine learning is embedded within optimization, the optimizer
does not use predictions uniformly: it preferentially explores
configurations that appear attractive according to the surrogate. As
a result, predictive performance on a conventional random holdout
sample can provide an incomplete picture of how the model will
perform when used to guide decisions. In the case study, prediction
and ranking performance deteriorate among high-scoring
search-generated candidates, while the component experiments show
meaningful differences in exact downstream savings across surrogate
and scoring choices. For managers, the implication is that AI
components used within optimization should be evaluated end to end,
using the quality of the resulting decisions in addition to
conventional predictive metrics.

\textbf{Insight 4: Learning and exact optimization play complementary roles
in strategic design.}
Learning enables the exploration of a much larger set of candidate
network structures than could practically be evaluated using the
exact lower-level model. Exact optimization can then be used to
recombine promising configurations, while exact lower-level
evaluation provides a reliable assessment of the selected network.
This separation suggests a practical architecture for AI-supported
strategic decisions: use learning where computational speed and broad
exploration are most valuable, while retaining optimization and exact
evaluation where feasibility and reliable performance assessment are critical.

The framework incorporates practical business requirements such as
cluster-size limits, distance restrictions, and central-site
eligibility rules, and can be extended to accommodate additional
considerations such as capacity restrictions, regional boundaries,
and service-time requirements. This flexibility is important because
implementable strategic network designs must satisfy operational and
organizational requirements in addition to economic objectives.

\subsection{Limitations and future research}
\label{subsec:future-research}

Several directions for future research follow from this study. First, the
lower-level evaluator can be extended to incorporate additional operational
realism, including implementation frictions associated with workload
concentration, warehousing complexity, transportation coordination, and
more dynamic inventory responses. Alternative lower-level models, including
richer shipment or expediting decisions and more flexible multi-echelon
structures, would also allow the robustness of the strategic design
recommendations to be assessed under different operating assumptions.

Second, the present computational study focuses on fulfillment centers.
Applying the framework to broader networks that also include middle-mile
and last-mile facilities would provide a richer setting in which to study
how strategic topology interacts with different operational roles and
constraints across the supply chain.

Third, the distribution-shift results suggest opportunities for more
adaptive learning procedures. The search can generate candidate structures
that differ from the configurations used to train the initial
surrogate, particularly when moving across geographic or operating
regimes. Future work could therefore combine surrogate-guided search with
selective exact evaluation of informative search-generated candidates,
followed by periodic model updating. Such an active-learning or
retraining mechanism could improve robustness while controlling the number
of additional exact lower-level evaluations.

Finally, the broader methodological framework extends beyond spare parts
applications. Many strategic operations problems require decisions over
large combinatorial design spaces while evaluating each candidate requires
solving an expensive operational optimization model. The combination of
graph-based representation, conservative surrogate scoring, search, and
exact final validation provides a general approach for such settings.

\subsection{Concluding remarks}

The main contribution of the paper is not to replace exact
optimization with machine learning, but to combine them in a way that
exploits their respective strengths. The graph-based surrogate
enables scalable exploration of candidate network structures,
conservative scoring limits the influence of optimistic predictions
during search, and optimization recombines promising configurations
into a feasible network. Exact lower-level evaluation then provides a
rigorous assessment of the final recommendation. This combination
makes large-scale exploration of strategic spare parts network
designs computationally practical while retaining exact evaluation
where it matters most. The Amazon case study demonstrates the
practical value of this approach and provides evidence that
learning-augmented optimization can improve strategic network design
while maintaining high service levels.

\newpage
\bibliographystyle{pomsref}

\let\oldbibliography\thebibliography
\renewcommand{\thebibliography}[1]{%
  \oldbibliography{#1}%
  \baselineskip14pt
  \setlength{\itemsep}{10pt}%
}
\bibliography{bibliography}

\newpage
\appendix
\input{arxiv_EC}

\end{document}

%% file: arxiv_EC.tex

\section{Surrogate Model Implementation Details}
\label{app:gnn-details}

This appendix provides the implementation
details of the graph
neural network (GNN) surrogate introduced in
Section~\ref{sec:solution-framework}. The purpose of the surrogate
is to approximate the lower-level evaluation model at negligible
computational cost during the search.

\subsection{Architecture and outputs}
\label{app:gnn-architecture}

We implement an attention-based graph neural network of the GATv2
family. Let $G=(V,E)$ denote a cluster graph, with node-feature
matrix $V^{feat} \in \mathbb{R}^{|V|\times d_v}$ and edge-feature
matrix $E^{\mathrm{feat}} \in \mathbb{R}^{|E|\times d_e}$. The
architecture consists of four graph-convolution layers with hidden
dimension 32. The first two layers use multi-head attention with
two heads, while the remaining two layers use a single head. After
each of the first two graph-convolution layers, we apply dropout
regularization with rate 0.2. Throughout the network we use the
SiLU activation function.

The hidden node representations evolve as
\[
  H^{(0)} = V^{feat},
\]
and
\[
  H^{(\ell+1)}
  =
  \mathrm{SiLU}\!\left(
    \mathrm{GATv2}^{(\ell)}\!\left(H^{(\ell)},E^{\mathrm{feat}}\right)
  \right),
  \qquad \ell=0,1,2,3.
\]

After the last graph-convolution layer, node embeddings are
aggregated through global add pooling:
\[
  h_G = \sum_{i \in V} H_i^{(4)}.
\]
The graph-level vector $h_G$ is then passed to a two-layer
multilayer perceptron with dimensions
\[
  32 \rightarrow 16 \rightarrow d_{\mathrm{out}},
\]
again with SiLU activations.

In the implementation used for the experiments, the surrogate
predicts four lower-level business metrics:
\[
  \hat y =
  \bigl(
    \widehat{po}^{\,\mathrm{n}},
    \widehat{stock}^{\,\mathrm{n}},
    \widehat{po}^{\,\mathrm{t}},
    \widehat{stock}^{\,\mathrm{t}}
  \bigr),
\]
where $n$ and $t$ indicate the metrics for the two-echelon and the
one-echelon configurations, respectively. Thus, $d_{\mathrm{out}} =
4$. Service level is not included among the surrogate targets
because it varies little across the evaluated cluster
configurations and remains consistently high in the lower-level solutions.

\subsection{Training details}
\label{app:gnn-training}

The full graph dataset contains 3{,}513 cluster graphs generated
from five subregions. We split this dataset into 80\% training
data, 10\% validation data, and 10\% test data. The split is
defined at the cluster level so as to prevent leakage across observations.

All node and edge features are standardized using a standard scaler
fitted on the training set. The model is trained with
mean-squared-error loss using AdamW with learning rate $10^{-3}$
and weight decay $10^{-4}$. We use mini-batches of size 32 and
gradient accumulation to obtain an effective batch size of 64.
Learning-rate scheduling is performed through OneCycleLR with
cosine annealing. Early stopping is applied with patience of 50
epochs based on validation-set performance.

To estimate predictive uncertainty, we train an ensemble of $K=5$
GNNs with independent random initialization on resampled versions
of the training set. For each member $k$, we convert the four predicted
cost components into the member-level net economic value
\[
  \widehat{\delta}^{(k)}(G)
  = \widehat{po}^{t,(k)}(G)+\widehat{so}^{t,(k)}(G)
  -\widehat{po}^{n,(k)}(G)-\widehat{so}^{n,(k)}(G).
\]
The conservative score used in the main text is then
\[
  \widehat{\delta}_{\tau}(G)
  =
  \mathrm{Quantile}_{\tau}
  \Bigl(
    \{\widehat{\delta}^{(k)}(G)\}_{k=1}^K
  \Bigr),
  \qquad \tau=0.3.
\]
Because larger $\delta$ indicates greater savings, this lower-quantile
aggregation penalizes candidate configurations whose attractive value
is not supported consistently across ensemble members.

\section{Candidate Library Generation, Initial Solution Completion,
Search, and Final Recombination}
\label{app:search-overview}

This appendix describes the generation of candidate
cluster library and the completion procedure used after the
initial set partitioning step., Next, it summarizes the surrogate-assisted
variable neighborhood search, and documents the final
set-partitioning recombination and exact validation steps.

\subsection{Generation of the candidate cluster library}
\label{ec:candidate-generation}

The initial candidate library is generated by repeatedly solving a
$p$-center clustering model over the eligible sites in each of five
subregions. For a specified number of clusters $p$, the model partitions
the eligible sites into $p$ clusters and minimizes the maximum distance
between a site and its assigned center, subject to minimum and maximum
cluster-size requirements. Sites excluded by the business rules are
removed before solving.

To generate structurally diverse candidates, we solve the model using
five alternative measures of proximity: common spare parts,
pairwise savings, non-common spare parts, geographic distance measured
using the Haversine formula, and geographic distance with an additional
upper bound on the maximum admissible central-to-leaf distance. We vary
both the number of clusters and the cluster-size bounds across multiple
configurations. Infeasible configurations are discarded. The resulting
pool is supplemented with randomly generated feasible partitions and
cluster samples.

Every cluster appearing in a feasible partition is extracted as an
individual candidate. We additionally generate central-site variants
by allowing each eligible leaf to serve as the central site while
preserving the cluster composition. Candidates with identical site
composition and central-site assignment are then removed.

Each remaining candidate is evaluated using the exact lower-level
model of Section~\ref{sec:nodal-model}. Candidates with invalid
evaluations or violating the final size, distance, or site-eligibility
requirements are discarded. After combining candidates generated
across the five subregions and removing duplicates, the resulting
library contains 3,513 unique feasible clusters. These clusters provide
the labeled data used to train and evaluate the surrogate and the
initial candidate pool used by the set-partitioning model.

\subsection{Completion of the initial set partitioning solution}
\label{app:initial-solution}

The set partitioning model of Section~\ref{subsec:set-partitioning}
is solved over a finite library of candidate clusters. Since this
library is not exhaustive, the optimal solution of the set
partitioning problem may cover only a subset of sites. The
procedure below is therefore used to increase coverage and
construct the initial solution passed to the metaheuristic.

Let $\mathcal{I}$ denote the full set of sites and let
$\mathcal{C}^\star \subseteq \mathcal{C}$ denote the subset of
clusters selected by the set partitioning model. Let
$\mathcal{I}^\star = \bigcup_{c \in \mathcal{C}^\star} c.\mathrm{sites}$
be the set of covered sites. The uncovered sites are then
$\hat{\mathcal{I}} = \mathcal{I} \setminus \mathcal{I}^\star$.

We first partition $\hat{\mathcal{I}}$ into two subsets:
\begin{align}
  \mathcal{I}_1
  &=
  \left\{
    i \in \hat{\mathcal{I}} \;:\;
    d_{ij} > d,\ \forall j \in \hat{\mathcal{I}}\setminus\{i\}
  \right\},
  \\
  \mathcal{I}_2
  &=
  \hat{\mathcal{I}} \setminus \mathcal{I}_1,
\end{align}
where $d$ is the maximum admissible distance between a central site
and a leaf site. Thus, $\mathcal{I}_1$ contains uncovered sites
that are isolated relative to the remaining uncovered sites, while
$\mathcal{I}_2$ contains the non-isolated uncovered sites. If
$|\mathcal{I}_2|<m$, where $m$ is the minimum cluster size, we set
$\mathcal{I}_1 \leftarrow \mathcal{I}_1 \cup \mathcal{I}_2$ and
$\mathcal{I}_2 \leftarrow \emptyset$.

Let $u=17$ denote the maximum admissible cluster size, let
$\mathcal{H}\subseteq\mathcal{I}$ denote the sites eligible to serve as
central sites after applying the business exclusions, and let $h(c)$
denote the selected center of cluster $c$.

If $\mathcal{I}_2 \neq \emptyset$, we solve an auxiliary covering
model on $\mathcal{I}_2$ in order to construct additional provisional
clusters. Let binary variable $y_i$ equal 1 if site $i$ is selected
as a cluster center, and let binary variable $x_{ij}$ equal 1 if
site $j$ is assigned to center $i$. We solve:
\begin{align}
  \min\quad & \sum_{i \in \mathcal{I}_2} y_i
  \label{eq:cover-obj}
  \\
  \text{s.t.}\quad
  & d_{ij}x_{ij} \le d\,y_i,
  \qquad \forall i,j \in \mathcal{I}_2,
  \label{eq:cover-2}
  \\
  & x_{ij} \le y_i,
  \qquad \forall i,j \in \mathcal{I}_2,
  \label{eq:cover-3}
  \\
  & x_{ii} \ge y_i,
  \qquad \forall i \in \mathcal{I}_2,
  \label{eq:cover-4}
  \\
  & \sum_{i \in \mathcal{I}_2} x_{ij} = 1,
  \qquad \forall j \in \mathcal{I}_2,
  \label{eq:cover-5}
  \\
  & x_{ij},y_i \in \{0,1\},
  \qquad \forall i,j \in \mathcal{I}_2.
  \label{eq:cover-6}
\end{align}

Let $\tilde{\mathcal{C}}_2$ denote the set of clusters induced by
the solution. We then define
\[
  \mathcal{C}_2 =
  \{ c \in \tilde{\mathcal{C}}_2 : m \le |c| \le u,
  \ h(c) \in \mathcal{H} \},
\]
and collect the sites belonging to clusters that fail the size or
center-eligibility requirements into
\[
  \mathcal{I}_3
  =
  \bigcup_{c \in \tilde{\mathcal{C}}_2 \setminus \mathcal{C}_2}
  c.\mathrm{sites}.
\]

The remaining uncovered sites are
\[
  \mathcal{I}_{\mathrm{remain}} = \mathcal{I}_1 \cup \mathcal{I}_3.
\]
For each site $i \in \mathcal{I}_{\mathrm{remain}}$, we attempt to
assign it to the nearest feasible central site among the selected
clusters in $\mathcal{C}^\star \cup \mathcal{C}_2$:
\[
  j^\star
  =
  \arg\min_{j \in H(\mathcal{C}^\star \cup \mathcal{C}_2)}
  \left\{
    d_{ij} : d_{ij}\le d,\ |c(j)|<u
  \right\},
\]
where $H(\mathcal{C}^\star \cup \mathcal{C}_2)$ denotes the set of
central sites of clusters in $\mathcal{C}^\star \cup
\mathcal{C}_2$, and $c(j)$ is the cluster centered at $j$. A proposed
attachment is retained only if it also satisfies the applicable
site-specific business restrictions. If no feasible assignment exists,
the site remains
unassigned and is handled later by the metaheuristic.

The resulting initial solution is therefore
\[
  S_0 = \mathcal{C}^\star \cup \mathcal{C}_2,
\]
possibly augmented with the reassigned sites from
$\mathcal{I}_{\mathrm{remain}}$.

\subsection{Variable neighborhood search: top-level algorithm}
\label{app:vns-main}

We use a surrogate-assisted variable neighborhood search (VNS) to
improve the initial network design. The search alternates between a
shaking phase, which perturbs the current solution, and a
local-search phase, which improves the perturbed solution. Whenever
an improvement is found, the search returns to the first
neighborhood; otherwise it moves to the next one.

\begin{algorithm}[H]
  \footnotesize
  \caption{Variable Neighborhood Search for Network Design}
  \label{alg:vns}
  \begin{algorithmic}[1]
    \Require Initial solution $S_0$, maximum neighborhood index
    $k_{\max}$, maximum number of outer iterations $M$, objective $f$
    \Ensure Best solution $S^\star$ found under the surrogate objective $f$
    \State $S^\star \gets S_0$
    \State $f^\star \gets f(S^\star)$
    \For{$t=1,\dots,M$}
    \State $k \gets 1$
    \While{$k \le k_{\max}$}
    \State $S' \gets \mathrm{Shaking}(S^\star,k)$
    \If{$S' = \mathrm{null}$}
    \State $k \gets k+1$
    \State \textbf{continue}
    \EndIf
    \State $S'' \gets \mathrm{LocalSearch}(S')$
    \State $f'' \gets f(S'')$
    \If{$f'' > f^\star$}
    \State $S^\star \gets S''$
    \State $f^\star \gets f''$
    \State $k \gets 1$
    \Else
    \State $k \gets k+1$
    \EndIf
    \EndWhile
    \EndFor
    \State \Return $S^\star$
  \end{algorithmic}
\end{algorithm}

\subsection{Neighborhood structures}
\label{app:vns-neighborhood-summary}

The search uses four neighborhood structures.

\begin{itemize}
  \item \textbf{$N_1$ (Swap).} A leaf site is removed from a donor
    cluster and used as the new central site of a receiver cluster.
    Any sites in the receiver that become infeasible under the new
    central are reassigned to other feasible clusters.
  \item \textbf{$N_2$ (Double swap).} The swap move of $N_1$ is
    applied twice consecutively.
  \item \textbf{$N_3$ (Break cluster).} A sufficiently large
    cluster is split by extracting a new cluster of minimum size
    around a newly selected central site.
  \item \textbf{$N_4$ (Merge clusters).} One cluster is dissolved
    and its sites are reassigned to other feasible clusters
    according to a randomized multi-criteria rule.
\end{itemize}

After each shaking move, local search re-optimizes the identity of
the central site within each cluster by checking whether one of the
leaf sites provides a better feasible center under the surrogate
objective. The descriptions above summarize the operator families;
implementation-specific move generation and repair details are not
reproduced here.

\subsection{Final set partitioning and exact validation}
\label{app:final-sp}

After each search run, we collect the distinct cluster compositions
retained from the promising network structures visited by that arm and
form an enriched candidate pool. We then solve the selective-coverage
set partitioning formulation of
Section~\ref{subsec:set-partitioning} over this pool. For a mean-scoring
arm, each cluster coefficient is the mean of the member-level predicted
economic values $\widehat{\delta}^{(k)}$; for a conservative-scoring
arm, it is the lower-quantile value $\widehat{\delta}_{0.3}$ defined in
Section~\ref{subsec:pessimistic-ensemble}. The set partitioning model is
solved to optimality for these arm-specific surrogate coefficients.

This combinatorial optimization step is exact for the stated
surrogate-valued set partitioning model; it does not make every
search-generated cluster coefficient an exact lower-level evaluation.
The network selected by the set partitioning model is subsequently
re-evaluated using the exact lower-level model of
Section~\ref{sec:nodal-model}. Only this exact re-evaluation is used to
report the final network-level business performance.

\subsection{Search Convergence and Dispersion}
\label{app:vns-convergence}

The component analysis uses six independent full-pipeline replications
per arm. Each replication trains a five-member surrogate ensemble,
executes 100 independent VNS restarts with 100 search steps per
restart, and consolidates the retained clusters through the final set
partitioning model to produce one network. The reported arm-level
means and standard deviations are computed across these six final
networks.

We examine the convergence behavior of the variable neighborhood search
to assess whether the improvements reported in the main paper depend on
a small number of favorable search trajectories. The production
configuration is initialized from the Net-baseline solution and uses the
GNN ensemble with conservative scoring. For the convergence diagnostic
shown below, we examine 100 independent restarts, each with a budget of
100 search steps. The figure is used only as descriptive convergence
evidence; the main text separately reports six-replication component
estimates.

\begin{figure}[t!]
  \centering
  \caption{Convergence and dispersion of the variable neighborhood search.
    Left: normalized surrogate-objective trajectories for the ten
    best-performing restarts, scaled to a best observed value of 1.00.
    Right: distribution of final improvements relative to the
  Net-baseline across 100 independent restarts.}
  \label{fig:metaheuristic_ec}
  \includegraphics[width=0.49\textwidth]{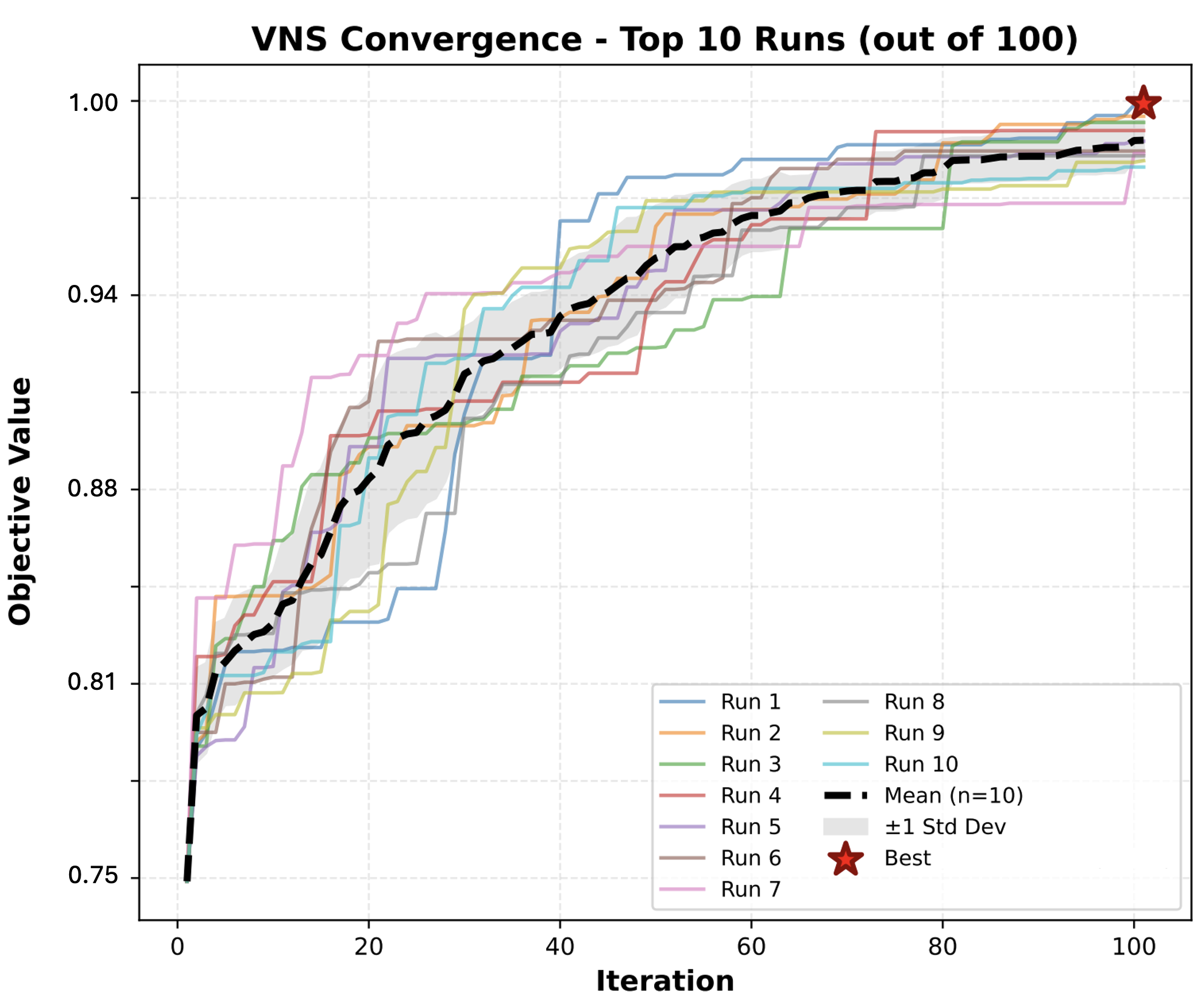}
  \includegraphics[width=0.47\textwidth]{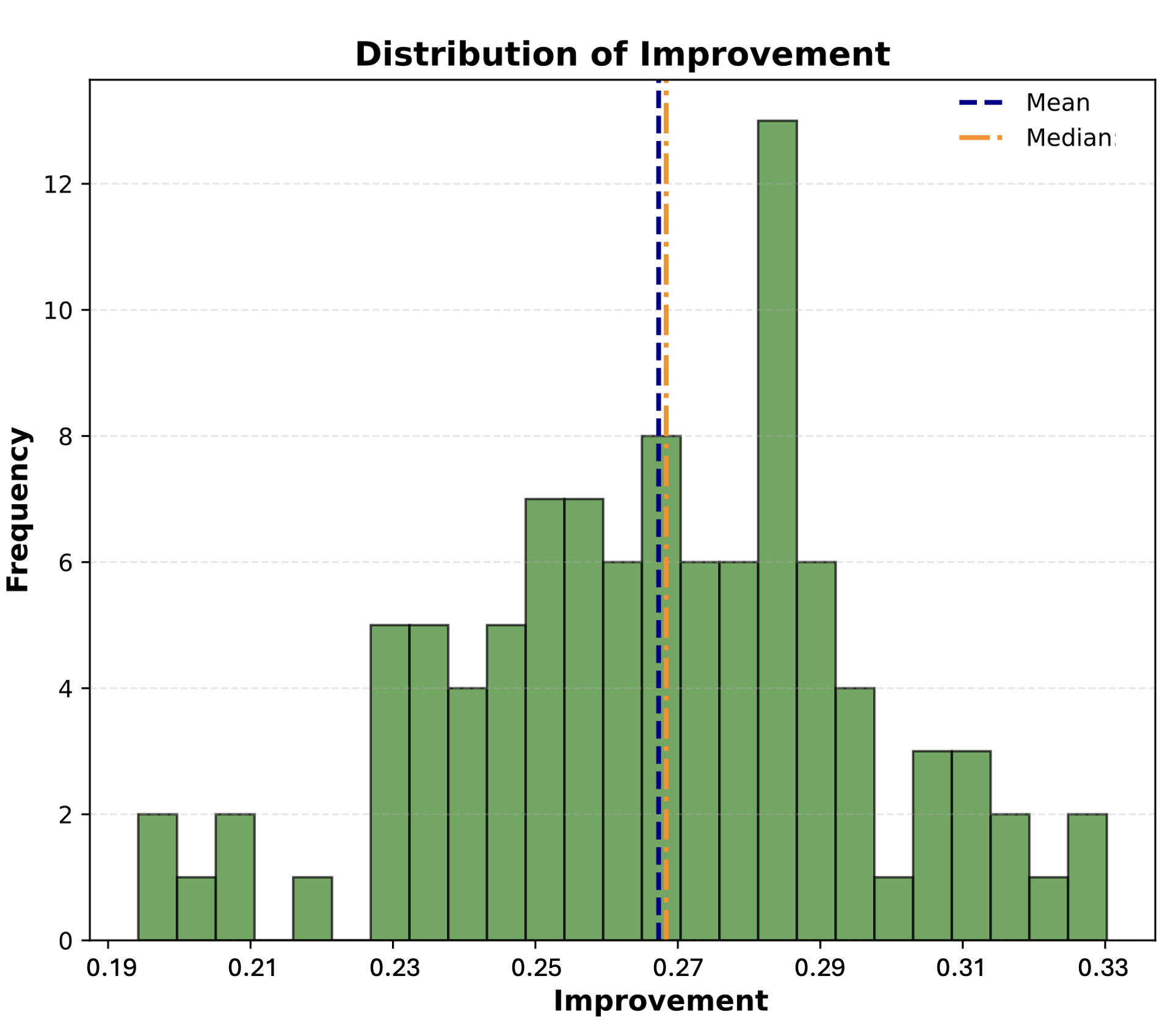}
\end{figure}

Figure~\ref{fig:metaheuristic_ec} summarizes the search trajectories.
The left panel reports the normalized surrogate-objective trajectories
for the ten best-performing restarts, scaled so that the best observed
objective value equals 1.00. The right panel shows the distribution of
final improvements relative to the initial Net-baseline across all 100
restarts. The selected top-ten trajectories illustrate convergence
behavior, whereas the all-restart histogram describes the dispersion
of final outcomes.

The top-ten trajectories show that high-quality solutions emerge
progressively as the search explores alternative cluster structures.
The all-restart histogram illustrates the dispersion induced by the
stochastic search. These plots provide descriptive evidence about
search behavior rather than a statistical stability guarantee.
This variability motivates the replicated experimental protocol used
in the component analysis and end-to-end evaluation in the main paper,
where reported network performance is based on independent search runs
and exact lower-level validation.

%% file: bibliography.bib
@misc{Battaglia2018,
  title = {Relational inductive biases, deep learning, and graph networks},
  author = {Battaglia, P. W. and Hamrick, J. B. and Bapst, V. and
            Sanchez-Gonzalez, A. and Zambaldi, V. and others},
  howpublished = {arXiv preprint arXiv:1806.01261},
  year = {2018},
}

@article{Cappart2021,
  title = {Combinatorial optimization and reasoning with graph neural networks},
  author = {Cappart, Q. and Ch{\'e}telat, D. and Khalil, E. B. and Lodi, A. and
            Morris, C. and Veli{\v{c}}kovi{\'c}, P.},
  journal = {Journal of Machine Learning Research},
  volume = {24},
  number = {130},
  pages = {1--61},
  year = {2023},
}

@article{kennedy_overview_2002,
  title = {An overview of recent literature on spare parts inventories},
  volume = {76},
  number = {2},
  journal = {International Journal of Production Economics},
  author = {Kennedy, W. J. and Patterson, J. W. and Fredendall, L. D.},
  month = mar,
  year = {2002},
  pages = {201--215},
}

@article{zhang_spare_2021,
  title = {Spare parts inventory management: A literature review},
  volume = {13},
  number = {5},
  journal = {Sustainability},
  publisher = {Multidisciplinary Digital Publishing Institute},
  author = {Zhang, S. and Huang, K. and Yuan, Y.},
  month = jan,
  year = {2021},
  pages = {2460},
}

@article{bacchetti_spare_2012,
  series = {Special {Issue} on {Forecasting} in {Management} {Science}},
  title = {Spare parts classification and demand forecasting for stock control:
           {Investigating} the gap between research and practice},
  volume = {40},
  number = {6},
  journal = {Omega},
  author = {Bacchetti, A. and Saccani, N.},
  year = {2012},
  pages = {722--737},
}

@book{axsater_inventory_2015,
  address = {Cham},
  series = {International {Series} in {Operations} {Research} \& {Management} {
            Science}},
  title = {Inventory {Control}},
  volume = {225},
  language = {en},
  publisher = {Springer International Publishing},
  author = {Axsäter, S.},
  year = {2015},
}

@article{zheng_finding_1991,
  title = {Finding optimal (s, {S}) policies is about as simple as evaluating a
           single policy},
  volume = {39},
  number = {4},
  journal = {Operations Research},
  publisher = {INFORMS},
  author = {Zheng, Y.-S. and Federgruen, A.},
  month = aug,
  year = {1991},
  pages = {654--665},
}

@article{sherbrooke_metric_1968,
  title = {{METRIC}: A multi-echelon technique for recoverable item control},
  volume = {16},
  number = {1},
  journal = {Operations Research},
  publisher = {INFORMS},
  author = {Sherbrooke, C. C.},
  month = feb,
  year = {1968},
  pages = {122--141},
}

@book{muckstadt2005analysis,
  title = {Analysis and Algorithms for Service Parts Supply Chains},
  author = {Muckstadt, J. A.},
  year = {2005},
  publisher = {Springer},
  address = {New York},
  series = {Springer Series in Operations Research and Financial Engineering},
}

@article{eppen_noteeffects_1979,
  title = {Note—effects of centralization on expected costs in a multi-location
           newsboy problem},
  volume = {25},
  number = {5},
  journal = {Management Science},
  publisher = {INFORMS},
  author = {Eppen, G. D.},
  month = may,
  year = {1979},
  pages = {498--501},
}

@article{paterson_inventory_2011,
  title = {Inventory models with lateral transshipments: {A} review},
  volume = {210},
  number = {2},
  journal = {European Journal of Operational Research},
  author = {Paterson, C. and Kiesmüller, G. and Teunter, R. and Glazebrook, K.},
  year = {2011},
  pages = {125--136},
}

@article{grahovac_sharing_2001,
  title = {Sharing and lateral transshipment of inventory in a supply chain with
           expensive low-demand items},
  volume = {47},
  number = {4},
  journal = {Management Science},
  publisher = {INFORMS},
  author = {Grahovac, J. and Chakravarty, A.},
  year = {2001},
  pages = {579--594},
}

@article{wong_multi-item_2006,
  series = {Feature {Cluster}: {Heuristic} and {Stochastic} {Methods} in {
            Optimization}},
  title = {Multi-item spare parts systems with lateral transshipments and
           waiting time constraints},
  volume = {171},
  number = {3},
  journal = {European Journal of Operational Research},
  author = {Wong, H. and van Houtum, G. J. and Cattrysse, D. and Oudheusden, D.
            Van},
  year = {2006},
  pages = {1071--1093},
}

@article{axsater_new_2003,
  title = {A new decision rule for lateral transshipments in inventory systems},
  volume = {49},
  number = {9},
  journal = {Management Science},
  publisher = {INFORMS},
  author = {Axsäter, S.},
  year = {2003},
  pages = {1168--1179},
}

@article{tapia-ubeda_modelling_2020,
  title = {Modelling and solving spare parts supply chain network design
           problems},
  volume = {58},
  number = {17},
  journal = {International Journal of Production Research},
  publisher = {Taylor \& Francis},
  author = {Tapia-Ubeda, F. J. and Miranda, P. A. and Roda, I. and Macchi, M.
            and Durán, O.},
  year = {2020},
  pages = {5299--5319},
}

@inproceedings{sinha2016optimistic,
  author = {Sinha, A. and Malo, P. and Deb, K.},
  title = {Solving Optimistic Bilevel Programs by Iteratively Approximating
           Lower Level Optimal Value Function},
  booktitle = {Proceedings of the IEEE Congress on Evolutionary Computation
               (CEC)},
  pages = {1877--1884},
  year = {2016},
  publisher = {IEEE},
}

@article{shen2003joint,
  author = {Shen, Z.-J. M. and Coullard, C. and Daskin, M. S.},
  title = {A Joint Location--Inventory Model},
  journal = {Transportation Science},
  volume = {37},
  number = {1},
  pages = {40--55},
  year = {2003},
}

@inproceedings{sinha2018kriging,
  author = {Sinha, A. and Bedi, S. and Deb, K.},
  title = {Bilevel Optimization Based on Kriging Approximations of Lower Level
           Optimal Value Function},
  booktitle = {Proceedings of the IEEE Congress on Evolutionary Computation
               (CEC)},
  pages = {1--8},
  year = {2018},
  publisher = {IEEE},
}

@article{elmachtoub2022smart,
  title = {Smart “Predict, then Optimize”},
  author = {Elmachtoub, A. and Grigas, P.},
  journal = {Management Science},
  volume = {68},
  number = {1},
  pages = {9--26},
  year = {2022},
}

@inproceedings{dumouchelle2024neurbilo,
  title = {{Neur2BiLO}: Neural bilevel optimization},
  author = {Dumouchelle, J. and Julien, E. and Kurtz, J. and Khalil, E. B.},
  booktitle = {Advances in Neural Information Processing Systems 37},
  pages = {86688--86719},
  year = {2024},
  doi = {10.52202/079017-2752},
}

@article{jin_surrogate-assisted_2011,
  title = {Surrogate-assisted evolutionary computation: {Recent} advances and
           future challenges},
  volume = {1},
  number = {2},
  journal = {Swarm and Evolutionary Computation},
  author = {Jin, Y.},
  month = jun,
  year = {2011},
  pages = {61--70},
}

@article{mladenovic1997variable,
  title = {Variable neighborhood search},
  author = {Mladenovi{\'c}, N. and Hansen, P.},
  journal = {Computers \& Operations Research},
  volume = {24},
  number = {11},
  pages = {1097--1100},
  year = {1997},
}

@article{hansen2010variable,
  title = {Variable neighborhood search: methods and applications},
  author = {Hansen, P. and Mladenovi{\'c}, N. and {Moreno P{\'e}rez}, J. A.},
  journal = {Annals of Operations Research},
  volume = {175},
  number = {1},
  pages = {367--407},
  year = {2010},
}

@article{Caserta2026,
  author = {Caserta, M. and D'Angelo, L.},
  title = {Optimising two-echelon spare parts inventory: A case study from {
           Amazon's} operations},
  year = {2026},
  pages = {1--19},
  journal = {International Journal of Production Research},
}

@article{ClarkScarf1960,
  author = {Clark, A. J. and Scarf, H.},
  title = {Optimal Policies for a Multi-Echelon Inventory Problem},
  journal = {Management Science},
  volume = {6},
  number = {4},
  pages = {475--490},
  year = {1960},
}

@article{Dada1992,
  author = {Dada, M.},
  title = {A Two-Echelon Inventory System with Priority Shipments},
  journal = {Management Science},
  volume = {38},
  number = {8},
  pages = {1140--1153},
  year = {1992},
}

@article{Driessen2015,
  author = {Driessen, M. A. and Arts, J. J. and van Houtum, G. J. and Rustenburg
            , W. D. and Huisman, B.},
  title = {Maintenance Spare Parts Planning and Control: A Framework for Control
           and Agenda for Future Research},
  journal = {Production Planning \& Control},
  volume = {26},
  number = {5},
  pages = {407--426},
  year = {2015},
}

@article{DrentArts2021,
  author = {Drent, M. and Arts, J.},
  title = {Expediting in Two-Echelon Spare Parts Inventory Systems},
  journal = {Manufacturing \& Service Operations Management},
  volume = {23},
  number = {6},
  pages = {1431--1448},
  year = {2021},
}

@article{Eruguz2016,
  author = {Eruguz, A. S. and Sahin, E. and Jemai, Z. and Dallery, Y.},
  title = {A comprehensive survey of guaranteed-service models for multi-echelon
           inventory optimization},
  journal = {International Journal of Production Economics},
  volume = {172},
  pages = {110--125},
  year = {2016},
}

@article{GravesWillems2000,
  author = {Graves, S. C. and Willems, S. P.},
  title = {Optimizing Strategic Safety Stock Placement in Supply Chains},
  journal = {Manufacturing \& Service Operations Management},
  volume = {2},
  number = {1},
  pages = {68--83},
  year = {2000},
}

@incollection{GravesWillems2003,
  author = {Graves, S. C. and Willems, S. P.},
  title = {Supply Chain Design: Safety Stock Placement and Supply Chain
           Configuration},
  booktitle = {Supply Chain Management: Design, Coordination and Operation},
  editor = {de Kok, A. G. and Graves, S. C.},
  series = {Handbooks in Operations Research and Management Science},
  volume = {11},
  pages = {95--132},
  publisher = {Elsevier},
  address = {Amsterdam},
  year = {2003},
}

@article{Howard2015,
  author = {Howard, C. and Marklund, J. and Tan, T. and Reijnen, I.},
  title = {Inventory Control in a Spare Parts Distribution System with Emergency
           Stocks and Pipeline Information},
  journal = {Manufacturing \& Service Operations Management},
  volume = {17},
  number = {2},
  pages = {142--156},
  year = {2015},
}

@article{Topan2010,
  author = {Topan, E. and Bay{\i}nd{\i}r, Z. P. and Tan, T.},
  title = {An Exact Solution Procedure for Multi-Item Two-Echelon Spare Parts
           Inventory Control Problem with Batch Ordering in the Central Warehouse
           },
  journal = {Operations Research Letters},
  volume = {38},
  number = {5},
  pages = {454--461},
  year = {2010},
}

@article{Topan2017,
  author = {Topan, E. and Bay{\i}nd{\i}r, Z. P. and Tan, T.},
  title = {Heuristics for Multi-Item Two-Echelon Spare Parts Inventory Control
           Subject to Aggregate and Individual Service Measures},
  journal = {European Journal of Operational Research},
  volume = {256},
  number = {1},
  pages = {126--138},
  year = {2017},
}

@article{Topan2020review,
  author = {Topan, E. and Eruguz, A. S. and Ma, W. and van der Heijden, M. C.
            and Dekker, R.},
  title = {A Review of Operational Spare Parts Service Logistics in Service
           Control Towers},
  journal = {European Journal of Operational Research},
  volume = {282},
  number = {2},
  pages = {401--414},
  year = {2020},
}

@article{TopanVanDerHeijden2020,
  author = {Topan, E. and van der Heijden, M. C.},
  title = {Operational Level Planning of a Multi-Item Two-Echelon Spare Parts
           Inventory System with Reactive and Proactive Interventions},
  journal = {European Journal of Operational Research},
  volume = {284},
  number = {1},
  pages = {164--175},
  year = {2020},
}

@article{Wang2019,
  author = {Wang, J. and Miao, H. and Yu, M.},
  title = {Interdependent Order Allocation in the Two-Echelon Competitive and
           Cooperative Supply Chain},
  journal = {International Journal of Production Research},
  volume = {57},
  number = {4},
  pages = {1190--1213},
  year = {2019},
}

@article{Wang2022,
  author = {Wang, Y. and Geunes, J. and Nie, X.},
  title = {Optimising Inventory Placement in a Two-Echelon Distribution System
           with Fulfillment-Time-Dependent Demand},
  journal = {International Journal of Production Research},
  volume = {60},
  number = {1},
  pages = {48--72},
  year = {2022},
}

@article{ChenDistribution2022,
  author = {Chen, H. and Dai, B. and Li, Y. and Zhang, Y. and Wang, X. and Deng,
            Y.},
  title = {Stock allocation in a two-echelon distribution system controlled by
           (s, {S}) policies},
  journal = {International Journal of Production Research},
  volume = {60},
  number = {3},
  pages = {894--911},
  year = {2022},
}

@article{Dui2023,
  author = {Dui, H. and Yang, X. and Liu, M.},
  title = {Importance Measure-Based Maintenance Analysis and Spare Parts Storage
           Configuration in Two-Echelon Maintenance and Supply Support System},
  journal = {International Journal of Production Research},
  volume = {61},
  number = {23},
  pages = {8325--8342},
  year = {2023},
}

@article{Geevers2024,
  author = {Geevers, K. and van Hezewijk, L. and Mes, M. R. K.},
  title = {Multi-Echelon Inventory Optimization Using Deep Reinforcement
           Learning},
  journal = {Central European Journal of Operations Research},
  volume = {32},
  number = {3},
  pages = {653--683},
  year = {2024},
}

@article{Stranieri2024,
  author = {Stranieri, F. and Stella, F. and Kouki, C.},
  title = {Performance of Deep Reinforcement Learning Algorithms in Two-Echelon
           Inventory Control Systems},
  journal = {International Journal of Production Research},
  volume = {62},
  number = {17},
  pages = {6211--6226},
  year = {2024},
}

@article{ChanLinSaxe2025,
  author = {Chan, T. C. Y. and Lin, B. and Saxe, S.},
  title = {Machine learning--augmented optimization of large bilevel and
           two-stage stochastic programs: Application to cycling network design},
  journal = {Manufacturing \& Service Operations Management},
  volume = {27},
  number = {6},
  pages = {1851--1868},
  year = {2025},
}

@article{SpieckermannMinnerSchiffer2025,
  author = {Spieckermann, C. and Minner, S. and Schiffer, M.},
  title = {Reduce-then-Optimize for the Fixed-Charge Transportation Problem},
  journal = {Transportation Science},
  volume = {59},
  number = {3},
  pages = {540--564},
  year = {2025},
}

@article{KraulSeizingerBrunner2023,
  author = {Kraul, S. and Seizinger, M. and Brunner, J. O.},
  title = {Machine Learning--Supported Prediction of Dual Variables for the
           Cutting Stock Problem with an Application in Stabilized Column
           Generation},
  journal = {INFORMS Journal on Computing},
  volume = {35},
  number = {3},
  pages = {692--709},
  year = {2023},
}

@article{caserta:dangelo:25,
  author = {Caserta, M. and D'Angelo, L.},
  title = {Intermittent demand forecasting for spare parts with little
           historical information},
  journal = {Journal of the Operational Research Society},
  volume = {76},
  number = {2},
  pages = {294--309},
  year = {2025},
}

@article{biggs+:23,
  author = {Biggs, M. and Hariss, R. and Perakis, G.},
  title = {Constrained optimization of objective functions determined from
           random forests},
  journal = {Production and Operations Management},
  volume = {32},
  number = {2},
  pages = {397-415},
  year = {2023},
}

@article{hammami:frein:14,
  author = {Hammami, R. and Frein, Y.},
  title = {A Capacitated Multi-echelon Inventory Placement Model under Lead Time
           Constraints},
  journal = {Production and Operations Management},
  volume = {23},
  number = {3},
  pages = {446-462},
  year = {2014},
}

@article{feng:shanthikumar:23,
  author = {Feng, Q. and Shanthikumar, J.},
  title = {The framework of parametric and nonparametric operational data
           analytics},
  journal = {Production and Operations Management},
  volume = {32},
  number = {9},
  pages = {2685-2703},
  year = {2023},
}
